\documentclass[11pt,letterpaper]{article}
\usepackage[margin=1in]{geometry}
\usepackage[T1]{fontenc}
\usepackage{lmodern}
\usepackage{amsmath,amssymb,amsfonts,amsthm,subfig,booktabs}
\numberwithin{equation}{section}
\usepackage{mathtools}
\usepackage[hidelinks]{hyperref}
\usepackage{etoolbox}
\usepackage{graphicx}
\usepackage{float}
\newtheorem{theorem}{Theorem}[section]
\newtheorem{corollary}{Corollary}[theorem]
\newtheorem{lemma}[theorem]{Lemma}
\newtheorem{remark}{Remark}

\newcommand{\dif}{\mathrm{d}}

\newcommand{\manuscriptfigure}[2]{\includegraphics[width=#1]{#2}}

\begin{document}

\date{}
\title{DMP-Oriented Anisotropic Mesh Adaptation for SUPG Finite Element Methods for Time-Dependent Convection-Diffusion Problems}
\author{Xianping Li%
\thanks{College of Integrative Sciences and Arts, Arizona State University, AZ, 85212, USA (\textit{xianping@asu.edu})}
}
\maketitle

\begin{abstract}
We study discrete bound preservation for fully discrete piecewise-linear streamline upwind Petrov--Galerkin (SUPG) 
finite element approximations of time-dependent convection--diffusion equations on anisotropic triangular meshes. 
Although the continuous homogeneous-source problem satisfies a parabolic maximum principle, SUPG approximations 
can exhibit nonphysical undershoots and overshoots. We derive sufficient nodal discrete maximum principle (DMP) conditions 
that account for the stabilized mass matrix, the combined effects of physical and streamline diffusion, time-step restrictions, 
and transfer between successive meshes. The analysis identifies a SUPG-effective diffusion tensor and motivates an 
anisotropic metric aligned with its inverse. We construct a combined metric that retains the effective-tensor geometry 
while introducing solution-dependent mesh density through recovered Hessian information. 
The resulting DMP conditions provide a theoretical basis for anisotropic mesh design. 
Numerical experiments for time-dependent transport problems with isotropic and moderately anisotropic diffusion demonstrate that,
under mass-compatible limited SUPG stabilization, the combined metric substantially reduces measured undershoots relative to other methods 
at matched nominal mesh budgets, although improved bound behavior can carry an accuracy penalty. 
A controlled tensor--density comparison in the anisotropic transport test supports a contribution 
from the effective tensor beyond the prescribed mesh density. 
\end{abstract}

\section{Introduction}
We consider the time-dependent convection--diffusion problem on a time-independent polygonal 
domain \(\Omega\subset\mathbb R^2\), with \(\Omega_T=\Omega\times[0,T]\):
\begin{subequations}
\label{modeleqn}
\begin{align}
\label{modelequation1}
&u_t-\varepsilon\nabla\cdot(\mathbb D\nabla u)+\boldsymbol b(\boldsymbol x,t)\cdot\nabla u=f(\boldsymbol x,t) &&\text{in }\Omega\times(0,T],\\
\label{Dirichletboundary}
&u(\boldsymbol x,t)=g(\boldsymbol x,t) &&\text{on }\partial\Omega\times[0,T],\\
\label{modeleqnboundary}
&u(\boldsymbol x,0)=u_0(\boldsymbol x) &&\text{in }\Omega.
\end{align}
\end{subequations}
Here \(\mathbb D=\mathbb D(\boldsymbol x)\) is symmetric and uniformly positive definite, 
with \(\lambda_{\min}(\mathbb D(\boldsymbol x))\ge d_0>0\) almost everywhere, 
\(\boldsymbol b\in L^\infty(\Omega_T)^2\), and \(\varepsilon\in(0,1]\) scales the physical diffusion. 
We assume sufficient regularity of the solution and data for the weak formulation and the chosen time discretization. 
No divergence-free condition is imposed unless explicitly stated.

For $f=0$, the parabolic comparison principle implies
\begin{equation}\label{eq:continuous-maximum-principle}
m_*\le u(\boldsymbol x,t)\le M_*
\quad\text{for a.e. }(\boldsymbol x,t)\in\Omega\times(0,T],
\end{equation}
provided the initial and Dirichlet data lie in $[m_*,M_*]$
\cite{protter1984maximum,lieberman1996parabolic}.
This holds for uniformly positive definite anisotropic diffusion
without requiring $\nabla\cdot\boldsymbol b=0$ in the nonconservative
convection form. Nonzero sources require source-dependent bounds.
A discrete maximum principle (DMP) requires the numerical solution,
under the corresponding assumptions, to preserve the bounds set by
admissible initial and boundary data. Standard Galerkin and linear SUPG
discretizations do not necessarily have this property
\cite{ciarlet1973maximum,li2013maximum,barrenechea2024dmp}.

When convection dominates diffusion, standard Galerkin finite element approximations 
can exhibit nonphysical oscillations unless the mesh is extremely fine. 
The resulting numerical difficulties are discussed systematically in \cite{roos2008robust}. 
SUPG stabilization adds dissipation preferentially in the streamline direction 
and remains one of the most widely used linear stabilization techniques for 
convection-dominated problems \cite{brooks1982streamline,codina2000stabilized,lube2006residual,cengizci2024supg}. 
Its oscillations near layers and the use of supplementary crosswind diffusion were studied by John and Knobloch \cite{john2007sold,john2008sold}. 
At the same time, anisotropic mesh adaptation can concentrate degrees of freedom near layers 
and align elements with strongly directional solution or coefficient structure \cite{nguyen2009adaptive,li2010anisotropic,Huang_book,apel1999anisotropic}. 
Time-dependent adaptive finite element methods are also well established \cite{eriksson1991adaptive,thomee2006parabolic}. 
Di, Xie, and Yin \cite{di2018anisotropic} developed a strategy that jointly designs anisotropic meshes 
and SUPG stabilization parameters for a steady two-dimensional convection--diffusion equation with 
scalar diffusivity and constant convection. Their metric and stabilization parameter are motivated 
by discretization-error control, rather than by a fully discrete DMP. Sun, Chen, and Xu \cite{sun2010numerical} 
showed through numerical studies of steady convection--diffusion--reaction problems that 
accurate layer resolution alone need not prevent oscillations, and that mesh structure and 
streamline stabilization both influence numerical stability.

In this work, we study conditions under which the SUPG finite element solution preserves nodal bounds. 
Discrete maximum principles are closely connected with monotonicity of the algebraic system and are strongly affected by mesh geometry. 
Classical finite element DMP results already show the central role of acute or non-obtuse simplicial meshes 
and matrix monotonicity \cite{ciarlet1973maximum,brandts2008dmp}. Xu and Zikatanov \cite{xu1999monotone} obtained an $M$-matrix structure 
for an edge-averaged convection--diffusion discretization under mesh assumptions. For anisotropic diffusion, 
Li and Huang derived an anisotropic non-obtuse-angle condition and corresponding metric tensors for DMP satisfaction 
and for the combination of DMP alignment with interpolation-error adaptation \cite{li2010anisotropic}. 
In two dimensions, weaker Delaunay-type conditions are also available for anisotropic diffusion \cite{huang2011delaunay}. 
Time-dependent anisotropic diffusion and diffusion--convection--reaction extensions were studied in \cite{li2013maximum,lu2014maximum}. 
Other approaches enforce discrete bounds by modifying the finite element
operator. Burman and Ern \cite{burman2002nonlinear} established a DMP for
nonlinear stabilized Galerkin methods on strictly acute meshes.
Kuzmin, Shashkov, and Svyatskiy \cite{kuzmin2009constrained} developed
a constrained method satisfying a DMP for anisotropic diffusion, while
Kuzmin \cite{kuzmin2006flux} proposed positivity-preserving algebraic
flux correction for scalar convection. These guarantees hold under the
respective assumptions of those methods. They do not establish a DMP
for the linear SUPG formulation considered in this paper. A recent survey
reviews the remaining challenges in combining accuracy and discrete
bound preservation \cite{barrenechea2024dmp}. Algebraic stabilization
provides further alternatives \cite{knobloch2023smuas}. Other constrained and invariant-domain-preserving formulations 
illustrate the distinct operator-based approach to enforcing bounds \cite{mudunuru2016maximum,guermond2017invariant}.

In contrast, we study DMP-oriented anisotropic mesh adaptation for a
specified linear SUPG discretization, without introducing an additional
nonlinear bound-enforcing correction to its operator. The theory provides
sufficient DMP conditions. The numerical experiments assess the reduction
of bound violations when those conditions are not fully met.
The principal geometric quantity is the effective tensor formed by
combining physical diffusion with streamline stabilization. We define
this tensor precisely from the element matrices in Section~\ref{subsec:dmp-preliminaries}
and use it to motivate anisotropic mesh metrics in Section~\ref{subsec-metric}.

Unlike the steady anisotropic SUPG metric/stabilization co-design of \cite{di2018anisotropic}, 
our analysis explicitly includes the SUPG time-derivative mass matrix, tensor-valued physical diffusion, 
time-step restrictions, and transfers between accepted meshes. We first derive explicit sufficient conditions 
under which the fully discrete linear SUPG scheme is monotone: an \(\mathbb H_K\)-metric angle condition controls 
the stiffness off-diagonal entries, a stabilization restriction controls the SUPG mass matrix, 
and time-step bounds yield the required \(M\)-matrix structure. These conditions provide a theoretical guarantee 
when all their hypotheses hold at each accepted time step, whereas practical anisotropic meshes 
typically fail one or more of the geometric and algebraic restrictions. 
The fact that a sufficient condition is violated gives no converse conclusion: 
undershoot may or may not occur. Accordingly, the numerical study distinguishes checking 
the theorem's hypotheses from comparing the observed magnitude of undershoot and overshoot across metric constructions, 
prioritizing full-trajectory violations over final-time extrema. We use mass-compatible limited stabilization 
for the principal transport comparisons. These comparisons evaluate coupled metric--stabilization strategies, 
not geometry in isolation. The effective-tensor metric and its solution-adaptive extension are
developed in Section~\ref{subsec-metric}. They retain the effective
diffusion geometry while permitting different element-density distributions.

The mesh conditions are independent of the mesh generator. We use BAMG
\cite{hecht2006bamg} for anisotropic adaptation and transfer old-time
solutions by $P_1$ interpolation, followed in the transport tests by a
sign-aware mass correction. The DMP theorem requires interval-preserving
transfer and does not follow from mass conservation alone. The need to account for modifications of 
finite element meshes over time has long been recognized \cite{dupont1982mesh}.

The rest of the paper is organized as follows. Section~\ref{sec:method}
describes the SUPG formulation, the adaptation-only metric, and solution
transfer. Section~\ref{sec:dmp-general} establishes sufficient DMP
conditions and derives the effective-tensor-based metrics.
Section~\ref{sec:numerical-experiments} presents the numerical results.
Conclusions are summarized in Section~\ref{sec:conclusions}.

\section{SUPG finite element discretization and anisotropic mesh adaptation}
\label{sec:method}
We use conforming, nondegenerate triangulations \(\mathcal T_h\) of a polygonal domain and continuous affine (\(P_1\)) basis functions. 
The formulas below are written for general spatial dimension \(d\), although our examples focus on two dimensional. 
We first describe the SUPG operator whose algebraic properties will be analyzed.

\subsection{Residual-based SUPG formulation}
Consider \eqref{modelequation1} with \eqref{Dirichletboundary}--\eqref{modeleqnboundary}. Let
\[
H_g^1(\Omega)=\{v\in H^1(\Omega):v|_{\partial\Omega}=g\},
\qquad
H_0^1(\Omega)=\{v\in H^1(\Omega):v|_{\partial\Omega}=0\}.
\]
The weak formulation of the problem is to find \(u\in H_g^1(\Omega)\) such that
\begin{equation}
\label{weakform}
(u_t,v)_\Omega+\varepsilon(\mathbb D\nabla u,\nabla v)_\Omega+(\boldsymbol b\cdot\nabla u,v)_\Omega=(f,v)_\Omega,
\qquad \forall v\in H_0^1(\Omega).
\end{equation}
For a fixed time \(t\), let \(V^h\subset H^1(\Omega)\) be the continuous piecewise-linear finite element space 
on \(\mathcal T_h\). We write \(U^h(t)=\{v_h\in V^h:v_h|_{\partial\Omega}=g_h(t)\}\), where \(g_h(t)\) is the nodal 
interpolant of the prescribed boundary data, and \(U_0^h=V^h\cap H_0^1(\Omega)\). We use the elementwise volume-averaged diffusion tensor
\begin{equation}
\label{volume-average-D}
\mathbb D_K=\frac1{|K|}\int_K\mathbb D(\boldsymbol x)\,d\boldsymbol x.
\end{equation}
Assume that $\mathbb D$ is measurable, symmetric, and uniformly elliptic, so that for some constants $0<d_0\le d_1<\infty$,
\begin{equation}
\label{ellipticity-D}
d_0|\boldsymbol\xi|^2\le\boldsymbol\xi^T\mathbb D(\boldsymbol x)\boldsymbol\xi
\le d_1|\boldsymbol\xi|^2,
\quad\boldsymbol\xi\in\mathbb R^d,\quad\text{for a.e. }\boldsymbol x\in\Omega.
\end{equation}
Integrating these inequalities over $K$ shows that $d_0I\preceq\mathbb D_K\preceq d_1I$. In particular, 
every element average is symmetric positive definite (SPD). Since the gradients of $P_1$ basis functions are constant on $K$,
\begin{equation}
\label{exact-averaged-p1}
\int_K(\nabla\varphi_i)^T\mathbb D(\boldsymbol x)\nabla\varphi_j\,d\boldsymbol x
=|K|(\nabla\varphi_i)^T\mathbb D_K\nabla\varphi_j.
\end{equation}
Thus, exact volume averaging does not alter the $P_1$ Galerkin diffusion stiffness matrix. 
When the element average is evaluated numerically, a positive-weight quadrature rule preserves SPD 
but introduces quadrature error unless it integrates $\mathbb D$ exactly. For \(w_h\in U^h\), define the element residual
\begin{equation}
\label{supg-residual}
R_K(w_h):=(w_h)_t-\varepsilon\nabla\cdot(\mathbb D_K\nabla w_h)+\boldsymbol b\cdot\nabla w_h-f.
\end{equation}
The residual-based SUPG approximation is defined by seeking \(u_h\in U^h\) such that
\begin{align}
\label{SUPGweakform}
&((u_h)_t,v_h)_\Omega
+\varepsilon\sum_{K\in\mathcal T_h}(\mathbb D_K\nabla u_h,\nabla v_h)_K
+(\boldsymbol b\cdot\nabla u_h,v_h)_\Omega \nonumber\\
&\quad +\sum_{K\in\mathcal T_h}\tau_K(R_K(u_h),\boldsymbol b\cdot\nabla v_h)_K
=(f,v_h)_\Omega,
\qquad \forall v_h\in U_0^h.
\end{align}
For \(P_1\) elements, \(\nabla u_h\) is constant on each element and therefore
\[
\nabla\cdot(\mathbb D_K\nabla u_h)=0
\qquad\text{in the interior of }K.
\]
Thus the element-interior diffusion part of the frozen-coefficient residual is zero. 
This fact does not make the frozen residual fully consistent for a variable physical tensor. 
Indeed, for a $P_1$ function $w_h$ and a sufficiently smooth tensor,
\[
-\varepsilon\nabla\cdot(\mathbb D\nabla w_h)
=-\varepsilon(\nabla\cdot\mathbb D)\cdot\nabla w_h
\quad\text{on }K,
\]
where $(\nabla\cdot\mathbb D)_j=\sum_i\partial_iD_{ij}$. Moreover, for a smooth exact solution $u$ 
of the physical PDE, the frozen strong residual satisfies
\begin{equation}
\label{frozen-residual-defect}
R_K(u)=\varepsilon\nabla\cdot[(\mathbb D-\mathbb D_K)\nabla u]
=\varepsilon\bigl((\mathbb D-\mathbb D_K):\nabla^2u
+(\nabla\cdot\mathbb D)\cdot\nabla u\bigr).
\end{equation}
Here the colon denotes the Frobenius double contraction,
$\mathbb A:\mathbb B=\sum_{i,j}A_{ij}B_{ij} =\operatorname{tr}(\mathbb A^T\mathbb B)$. 
The coefficient-gradient term in \eqref{frozen-residual-defect} need not vanish as $h\to0$, 
even though $\mathbb D_K$ approximates $\mathbb D$ locally. Hence exactness of the Galerkin diffusion matrix 
in \eqref{exact-averaged-p1} must not be confused with strong-residual consistency. 

Define
\begin{equation}
\label{Pe-tensor}
\mathrm{Pe}_K:=\frac{B_K h_K}{2\varepsilon\lambda_{\min}(\mathbb D_K)},
\qquad
B_K:=\|\boldsymbol b\|_{\infty,K},
\qquad h_K:=\operatorname{diam}(K),
\end{equation}
and use
\begin{equation}
\label{taudef}
\tau_K=\frac{h_K}{2B_K}\,\xi(\mathrm{Pe}_K),
\qquad
\xi(\mathrm{Pe}_K)=\min\left\{1,\frac{\mathrm{Pe}_K}{3}\right\},
\end{equation}
with \(\tau_K=0\) if \(B_K=0\). This nominal choice uses the element diameter and the smallest eigenvalue of 
the element-averaged physical diffusion tensor. For anisotropic meshes, other constructions are possible. 
In particular, Di, Xie, and Yin \cite[Eq.~(4.8)]{di2018anisotropic} proposed a stabilization parameter 
involving element dimensions along and across the convection direction, derived together with an 
error-oriented anisotropic metric for a steady scalar-diffusion problem. We do not identify their parameter 
with \(\tau_K\) in \eqref{taudef}, or assert that it satisfies our SUPG mass-sign restriction. 
In particular, the two branches of \eqref{taudef} have different scalings:
\[
\tau_K=\frac{h_K^2}{12\varepsilon\lambda_{\min}(\mathbb D_K)}
\quad(\mathrm{Pe}_K\le3),
\qquad
\tau_K=\frac{h_K}{2B_K}\quad(\mathrm{Pe}_K\ge3).
\]
We use \eqref{taudef} consistently as the nominal SUPG stabilization in the numerical comparisons. 
This scaling depends on element geometry, velocity, and diffusion, but not on the time-step size. 
No optimality claim is made for transient problems. A time-step-dependent alternative would also change 
the mass and stiffness conditions and is not a consequence of the present DMP theorem.
The sufficient DMP conditions developed in Section~\ref{sec:dmp-general} are stated directly 
in terms of a nonnegative \(\tau_K\) and can therefore be checked for other stabilization choices. 
Within each nominal or mass-limited numerical comparison, the same stabilization policy is applied across mesh metrics. 
Comparisons between policies also alter the mesh-dependent stabilized operators.

The resulting semidiscrete equations can be written in matrix form as follows, 
which makes the monotonicity conditions at each time step explicit.

\subsection{Matrix form and time discretization}
Write \(u_h(\boldsymbol x,t)=\sum_{j=1}^{N_v}u_j(t)\varphi_j(\boldsymbol x)\). 
Then the semi-discrete system is
\begin{equation}
\label{TimeODE}
\boldsymbol M\boldsymbol u_h'(t)+\boldsymbol A(t)\boldsymbol u_h(t)=\boldsymbol f(t),
\end{equation}
with
\begin{align}
\label{MmatrixFormula}
M_{ij}
&=\sum_{K\in\mathcal T_h}\left[(\varphi_j,\varphi_i)_K+\tau_K(\varphi_j,\boldsymbol b\cdot\nabla\varphi_i)_K\right],\\
\label{Amatrixint}
A_{ij}
&=\sum_{K\in\mathcal T_h}\left[(\boldsymbol b\cdot\nabla\varphi_j,\varphi_i)_K
+\varepsilon|K|(\nabla\varphi_i)^T\mathbb D_K\nabla\varphi_j
+\tau_K(\boldsymbol b\cdot\nabla\varphi_j,\boldsymbol b\cdot\nabla\varphi_i)_K\right],\\
\label{Frightside}
f_i(t)
&=\sum_{K\in\mathcal T_h}\left[(f,\varphi_i)_K+\tau_K(f,\boldsymbol b\cdot\nabla\varphi_i)_K\right].
\end{align}
The coefficient-dependent matrices are evaluated at
\[
t_{m+\theta}=(1-\theta)t_m+\theta t_{m+1}.
\]
Suppressing the stage superscript, the frozen-stage \(\theta\)-scheme is
\begin{align}
\label{discretizedtime}
\boldsymbol M\frac{\boldsymbol u^{m+1}-\boldsymbol u^m}{\Delta t}
+\boldsymbol A\bigl(\theta\boldsymbol u^{m+1}+(1-\theta)\boldsymbol u^m\bigr)
=\theta\boldsymbol f^{m+1}+(1-\theta)\boldsymbol f^m,
\end{align}
which is equivalent to
\begin{equation}
\label{ODEsystem}
\widetilde{\boldsymbol A}\boldsymbol u^{m+1}=\widetilde{\boldsymbol f}^{m+1},
\end{equation}
where
\begin{align}
\label{Amatrix}
\widetilde{\boldsymbol A}&=\boldsymbol M+\theta\Delta t\boldsymbol A,\\
\label{fvector}
\widetilde{\boldsymbol f}^{m+1}
&=\Delta t[\theta\boldsymbol f^{m+1}+(1-\theta)\boldsymbol f^m]
+[\boldsymbol M-(1-\theta)\Delta t\boldsymbol A]\boldsymbol u^m.
\end{align}
The same stage-frozen values of \(\boldsymbol b\), \(\tau_K\), and the element tensor are used 
in the matrices on both sides of \eqref{discretizedtime}. The two load vectors in that equation are 
assembled at \(t_m\) and \(t_{m+1}\) using these stage-frozen streamline test weights. 
In particular, the coefficient matrices are evaluated at the midpoint when \(\theta=1/2\). 
For time-dependent velocity, this describes a particular stage-frozen time integrator. 
Temporal accuracy must be checked for this integrator rather than assumed from the classical autonomous \(\theta\)-scheme.

We next introduce the Hessian-based adaptation metric used as the accuracy-oriented baseline in the numerical comparisons.

\subsection{Hessian recovery and the adaptation-only metric}
\label{subsec:hessian-adapt}
We first consider a metric based solely on the recovered Hessian. Following the construction in 
Eq.~(68) of Li and Huang \cite{li2010anisotropic}, which is obtained by minimizing an interpolation-error bound 
in the \(H^1\)-seminorm, we define
\begin{equation}
\label{adapt-AK}
\mathbb A_K:=\mathbb I+\frac{1}{\alpha_h}|H_K(u_h)|
\end{equation}
and, in spatial dimension \(d\),
\begin{equation}
\label{adapt-rhoK}
\rho_K
:=\|\mathbb A_K\|_2^{d/(d+2)}
\det(\mathbb A_K)^{1/(d+2)}.
\end{equation}
Then the adaptation-only metric is
\begin{equation}
\label{metrictensordef}
\mathbb{M}_{\mathrm{adapt},K}
=\rho_K^{2/d}\det(\mathbb A_K)^{-1/d}\mathbb A_K .
\end{equation}
For the two-dimensional problems considered in this paper, it becomes
\begin{equation}
\label{metrictensordef-2d}
\mathbb{M}_{\mathrm{adapt},K}
=\|\mathbb A_K\|_2^{1/2}\det(\mathbb A_K)^{-1/4}\mathbb A_K.
\end{equation}
Notice that Eq.~(68) of \cite{li2010anisotropic} uses the Frobenius norm. We replace it by the 
spectral norm \(\|\cdot\|_2\) in \eqref{adapt-rhoK}. The positive regularization parameter \(\alpha_h\) 
is determined implicitly by the corresponding Li--Huang density normalization
\begin{equation}
\label{alpha-adapt-normalization}
\sum_{K\in\mathcal T_h}|K|\rho_K=2|\Omega|.
\end{equation}
If the recovered Hessian vanishes on all elements, the normalization equation \eqref{alpha-adapt-normalization} 
has no positive solution. In that degenerate case, we define \(\mathbb{M}_{\mathrm{adapt},K}=\mathbb I\) 
(up to a global complexity scaling) and do not solve the normalization equation. Otherwise, the existence and uniqueness 
of the normalization parameter follow directly from the eigenvalues of the regularized Hessian. 
Indeed, if $\lambda_{K,j}\ge0$ are the eigenvalues of $|H_K(u_h)|$, then
$\mathbb A_K$ has eigenvalues $1+\lambda_{K,j}/\alpha_h$.
As $\alpha_h$ increases, each eigenvalue is nonincreasing and at
least one is strictly decreasing on an element with nonzero Hessian.
Both $\|\mathbb A_K\|_2$ and $\det(\mathbb A_K)$ are positive and
nonincreasing, so $\rho_K$ is continuous and the total density is
strictly decreasing. Its limit is $|\Omega|$ as
$\alpha_h\to\infty$ and is $+\infty$ as $\alpha_h\downarrow0$.
The intermediate value theorem consequently gives exactly one
$\alpha_h>0$ satisfying \eqref{alpha-adapt-normalization}. Here \(H_K(u_h)\) is the recovered Hessian of 
the finite element solution and \(|H_K(u_h)|\) is its matrix absolute value, obtained by replacing 
the Hessian eigenvalues by their absolute values. The recovered Hessian is obtained by 
local quadratic least-squares fitting to the nodal values. Thus the metric retains the structure 
of Eq.~(68) in \cite{li2010anisotropic}, but differs in the choice of matrix norm and Hessian recovery.

Given a metric tensor, we construct a new mesh and transfer the accepted solution before advancing the PDE.

\subsection{BAMG anisotropic remeshing and solution transfer}
\label{subsec:bamg}
We generate anisotropic triangular meshes from an SPD metric field
using BAMG \cite{hecht2006bamg}, following the metric-based approach
of \cite{li2010anisotropic}. The mesh generator controls geometry and
complexity but does not enforce the algebraic DMP conditions.

Let \(\mathcal T_h^m\) denote the accepted mesh at time \(t_m\). When a new mesh \(\widehat{\mathcal T}_h\) is generated, 
the old-time finite element function is transferred by pointwise \(P_1\) interpolation. 
If \(P\) denotes the corresponding interpolation matrix, then each row consists of barycentric coordinates in 
the old triangle containing the new node. Consequently, for valid containing triangles,
\[
P\ge0,
\qquad
P\boldsymbol 1=\boldsymbol 1.
\]
Thus the interpolation preserves nodal bounds. It assumes that the new nodes lie in the closure of the old domain, 
the containing triangles are nondegenerate, and no extrapolation or nonpositive projection is used. 
It says nothing about conservation of integral mass or the interpolation error.

For the principal numerical transport comparisons, let
$\widehat{\boldsymbol u}=P\boldsymbol u^m$ be the barycentrically
interpolated nodal vector on the destination triangulation. Its piecewise-affine
finite element function is $\widehat u_h=\sum_i\widehat u_i\widehat\varphi_i$,
where $\widehat\varphi_i$ denotes the destination-mesh nodal basis. Let
$w_i=\int_\Omega\widehat\varphi_i\,\dif\boldsymbol x$ be its positive
$P_1$ mass weights. A sign-aware correction constructs
$\boldsymbol z$ satisfying
\begin{equation}\label{eq:mass-transfer-numerical}
 \sum_i w_i z_i=\int_\Omega u_h^m\,\dif\boldsymbol x,\qquad
 z_i=\widehat u_i\quad(i\in\mathcal V_{\mathrm{fix}}).
\end{equation}
Here $\mathcal V_{\mathrm{fix}}$ consists of boundary vertices and all
vertices of destination triangles where $\widehat u_h$ violates a
physical bound. On eligible interior vertices, a weighted quadratic correction minimizes
changes from the interpolated values subject to signed-mass conservation
and nodal bounds:
\begin{equation}\label{eq:transfer-objective}
\min_{\boldsymbol z}\frac12\sum_{i\in S}\frac{w_i}{m_i}
(z_i-\widehat u_i)^2,
\qquad \sum_i w_i z_i=\int_\Omega u_h^m\,\dif\boldsymbol x.
\end{equation}
Here $S$ is an eligible support, $m_i>0$ weights the correction, and other nodal values are fixed. 
The correction preserves global signed mass, not elementwise mass.
Its feasibility and bound preservation are not consequences of the
mass constraint alone. The nonlinear map meets the transfer hypothesis
of the DMP theorem only if it preserves admissible intervals.
For a different conservative interpolation construction see
\cite{alauzet2010conservative}.

\section{Conditions for discrete maximum principle satisfaction}
\label{sec:dmp-general}

In this section, we develop verifiable sufficient conditions on the mesh, stabilization parameter, 
and time step size under which the fully discrete SUPG approximation satisfies a discrete maximum principle (DMP). 
We first give deliberately conservative elementwise conditions and then an alternative sign criterion 
for the fully assembled matrices. It is known that standard linear finite element approximations of 
convection-dominated problems need not satisfy a DMP. See, e.g., \cite{barrenechea2024dmp,knobloch2023smuas}. 
We consider monotonicity of the algebraic system obtained at each time step.

We assume that the mesh has been generated and the previous solution has been transferred to it before the time step is taken. 
The finite element basis functions are therefore fixed during the step.
We begin with element identities and identify the tensor governing the diffusion and streamline parts of the SUPG stiffness matrix.

\subsection{Preliminaries and the effective SUPG diffusion tensor}
\label{subsec:dmp-preliminaries}

All elementwise quantities in this section are evaluated at a single frozen time stage. 
The assembled matrices are those of \eqref{MmatrixFormula}--\eqref{Amatrixint} before imposing Dirichlet equations. 
The local estimates apply to rows associated with interior test nodes and to columns associated with 
either interior or boundary trial nodes. For each simplex \(K\in\mathcal T_h\), let
\[
B_K:=\|\boldsymbol b\|_{\infty,K},
\qquad
\boldsymbol q_i^K:=\nabla\varphi_i|_K,
\qquad
|\boldsymbol q_i^K|=\frac1{h_i^K},
\]
where \(h_i^K\) is the altitude from vertex \(i\) to its opposite face. For \(P_1\) basis functions,
\begin{equation}
\label{basis-integrals-dmp}
\int_K\varphi_i\,d\boldsymbol x=\frac{|K|}{d+1},
\qquad
\int_K\varphi_i\varphi_j\,d\boldsymbol x=
\begin{cases}
\displaystyle\frac{|K|}{(d+1)(d+2)},&i\ne j,\\[1ex]
\displaystyle\frac{2|K|}{(d+1)(d+2)},&i=j.
\end{cases}
\end{equation}
Moreover,
\begin{equation}
\label{convection-inequality}
\left|\int_K\varphi_i(\boldsymbol b\cdot\boldsymbol q_j^K)\,d\boldsymbol x\right|
\le \frac{|K|}{d+1}B_K|\boldsymbol q_j^K|.
\end{equation}

Define the element average
\begin{equation}
\label{bb-average}
\overline{\boldsymbol b\boldsymbol b^T}_K
:=\frac1{|K|}\int_K\boldsymbol b(\boldsymbol x,t)\boldsymbol b(\boldsymbol x,t)^T\,d\boldsymbol x
\end{equation}
and the effective SUPG diffusion tensor
\begin{equation}
\label{H-eff}
\mathbb H_K:=\varepsilon\mathbb D_K+\tau_K\overline{\boldsymbol b\boldsymbol b^T}_K.
\end{equation}
Because \(\mathbb D_K\) is SPD and \(\overline{\boldsymbol b\boldsymbol b^T}_K\) is symmetric positive semidefinite,
\begin{equation}
\label{H-SPD}
\mathbb H_K\ \hbox{is SPD for every }\varepsilon>0\hbox{ and }\tau_K\ge0.
\end{equation}
Notice that this property does not require \(\nabla\!\cdot\boldsymbol b=0\).

Since the basis gradients are constant on \(K\), the local diffusion plus streamline contribution in \eqref{Amatrixint} is exactly
\begin{align}
&\varepsilon |K|(\boldsymbol q_i^K)^T\mathbb D_K\boldsymbol q_j^K
+\tau_K\int_K(\boldsymbol b\cdot\boldsymbol q_j^K)(\boldsymbol b\cdot\boldsymbol q_i^K)\,d\boldsymbol x
\nonumber\\
&\hspace{35mm}=|K|(\boldsymbol q_i^K)^T\mathbb H_K\boldsymbol q_j^K.
\label{H-exact}
\end{align}
For \(i\ne j\), define the metric angle \(\alpha_{ij}^{H,K}\) by
\begin{equation}
\label{H-angle}
\cos\alpha_{ij}^{H,K}
:=-\frac{(\boldsymbol q_i^K)^T\mathbb H_K\boldsymbol q_j^K}
{\|\mathbb H_K^{1/2}\boldsymbol q_i^K\|\,\|\mathbb H_K^{1/2}\boldsymbol q_j^K\|}.
\end{equation}
Thus
\begin{equation}
\label{H-cross}
(\boldsymbol q_i^K)^T\mathbb H_K\boldsymbol q_j^K
=-\|\mathbb H_K^{1/2}\boldsymbol q_i^K\|\,
 \|\mathbb H_K^{1/2}\boldsymbol q_j^K\|\cos\alpha_{ij}^{H,K}.
\end{equation}

We now derive conservative conditions ensuring that the stabilized mass entries are 
nonnegative and the stiffness off-diagonal entries are nonpositive.

\subsection{Spatial matrix conditions}
\label{subsec:spatial-dmp}

The first lemma bounds the streamline contribution to the consistent mass matrix.

\begin{lemma}[Nonnegativity and positive row sums of the SUPG mass matrix]
\label{lemma-mij}
Assume that, for every element $K$ and every local vertex $i$,
\begin{equation}
\label{cond-mij}
\tau_KB_K|\boldsymbol q_i^K|\le\frac1{d+2},
\qquad\text{equivalently}\qquad
\tau_KB_K\le\frac{h_i^K}{d+2}.
\end{equation}
Then all entries $M_{ij}$ with an interior test node $i$ and an
interior or boundary trial node $j$ are nonnegative. Moreover,
$M_{ii}>0$ and $\sum_j M_{ij}>0$ for every interior test node.
\end{lemma}

\begin{proof}
Let $M_{ij}^K$ denote the contribution from an element containing both
local vertices $i$ and $j$. Since $\boldsymbol q_i^K$ is constant, and by the velocity bound
$|\boldsymbol b|\le B_K$ and the exact $P_1$ identity
$\int_K\varphi_j\,d\boldsymbol x=|K|/(d+1)$, the mass formula gives
\begin{align*}
M_{ij}^K
&=\int_K\varphi_i\varphi_j\,d\boldsymbol x
  +\tau_K\boldsymbol q_i^K\cdot\int_K\varphi_j\boldsymbol b\,d\boldsymbol x,\\
\left|\tau_K\boldsymbol q_i^K\cdot
  \int_K\varphi_j\boldsymbol b\,d\boldsymbol x\right|
&\le\frac{|K|\tau_KB_K|\boldsymbol q_i^K|}{d+1}.
\end{align*}
For distinct vertices, the exact $P_1$ mass integral and
\eqref{cond-mij} therefore imply
\[
M_{ij}^K\ge\frac{|K|}{d+1}
\left(\frac1{d+2}-\tau_KB_K|\boldsymbol q_i^K|\right)\ge0.
\]
For the diagonal entry, the Galerkin mass contribution is twice as
large. Hence
\[
M_{ii}^K\ge\frac{|K|}{d+1}
\left(\frac2{d+2}-\tau_KB_K|\boldsymbol q_i^K|\right)
\ge\frac{|K|}{(d+1)(d+2)}>0.
\]
Assembly sums these nonnegative local contributions. A pair of nodes
with disjoint supports gives a zero entry. In particular, the result
includes columns belonging to Dirichlet boundary nodes, which must
not be discarded when the matrix signs are tested.

To prove strict positivity of the full row sum, use the partition
of unity on every element, including boundary shape functions:
\begin{align*}
\sum_jM_{ij}
&=\sum_{K\in\omega_i}\left(
\int_K\varphi_i\,d\boldsymbol x
+\tau_K\boldsymbol q_i^K\cdot\int_K\boldsymbol b\,d\boldsymbol x\right)\\
&\ge\sum_{K\in\omega_i}|K|
\left(\frac1{d+1}-\tau_KB_K|\boldsymbol q_i^K|\right)\\
&\ge\sum_{K\in\omega_i}\frac{|K|}{(d+1)(d+2)}>0.
\end{align*}
The final inequality uses $d+1<d+2$ and holds even if equality occurs
in some of the local restrictions \eqref{cond-mij}.
\end{proof}

The mass restriction depends on the smallest element altitude, while the
nominal stabilization parameter uses the diameter. 
A convenient choice that enforces the mass restriction is the
following limiter:
\begin{equation}
\label{tau-dmp-limiter}
\tau_K^{\mathrm DMP}:=
\min\left\{\tau_K^{\mathrm nominal},\ 
\min_{i\in\mathcal V(K)}\frac{h_i^K}{(d+2)B_K}\right\},
\end{equation}
with no restriction when $B_K=0$. This guarantees the needed mass-entry
signs but is not designed to optimize approximation accuracy.

In two dimensions, when $B_K>0$ and $\mathrm{Pe}_K\ge3$, the limiter gives
$\tau_K^{\mathrm{DMP}}/\tau_K^{\mathrm{nominal}}=\min\{1,h_{\min,K}/(2h_K)\}$. 
Thus highly stretched elements may require a substantial reduction of streamline stabilization. 
This restriction concerns the consistent SUPG mass matrix and is not a
limitation shared by every choice of mass treatment.

For \(i\ne j\), define
\begin{equation}
\label{delta-H}
\delta_{ij}^{K,H}
:=\|\mathbb H_K^{1/2}\boldsymbol q_i^K\|\,
  \|\mathbb H_K^{1/2}\boldsymbol q_j^K\|\cos\alpha_{ij}^{H,K}
-\frac{B_K}{d+1}|\boldsymbol q_j^K|.
\end{equation}

Then we have the following lemma.

\begin{lemma}[Off-diagonal stiffness estimate]
\label{lemma-aij}
If
\begin{equation}
\label{delta-H-cond}
\delta_{ij}^{K,H}\ge0
\qquad\text{for every }K\text{ and every ordered local pair }i\ne j,
\end{equation}
then $A_{ij}\le0$ for every distinct global pair of trial and test
nodes $i\ne j$.
\end{lemma}

\begin{proof}
For distinct local vertices, the $P_1$ basis gradients are constant.
By \eqref{H-exact}, the local stiffness entry is
\[
A_{ij}^K
=\int_K\varphi_i(\boldsymbol b\cdot\boldsymbol q_j^K)
\,d\boldsymbol x
+|K|(\boldsymbol q_i^K)^T\mathbb H_K\boldsymbol q_j^K.
\]
The first term can have either sign, but \eqref{convection-inequality}
gives the upper bound $|K|B_K|\boldsymbol q_j^K|/(d+1)$.
Using \eqref{H-cross} in the second term and then the definition of
$\delta_{ij}^{K,H}$, we obtain
\begin{align*}
A_{ij}^K
&\le |K|\left[\frac{B_K|\boldsymbol q_j^K|}{d+1}
-\|\mathbb H_K^{1/2}\boldsymbol q_i^K\|
 \|\mathbb H_K^{1/2}\boldsymbol q_j^K\|
\cos\alpha_{ij}^{H,K}\right]\\
&=-|K|\delta_{ij}^{K,H}\le0.
\end{align*}
The convection estimate differentiates the trial function
$\varphi_j$, which explains why the hypothesis is imposed on
ordered pairs and cannot be replaced by a symmetric condition
without further justification. Summing the local bounds over the
common element patch proves $A_{ij}\le0$. Nonoverlapping basis
functions give $A_{ij}=0$.
\end{proof}

Equivalently, the sufficient local angle condition is
\begin{equation}
\label{effective-angle-condition}
\cos\alpha_{ij}^{H,K}
\ge
\frac{B_K|\boldsymbol q_j^K|}
{(d+1)\|\mathbb H_K^{1/2}\boldsymbol q_i^K\|\,
       \|\mathbb H_K^{1/2}\boldsymbol q_j^K\|}.
\end{equation}
Notice that this condition is stronger than the nonobtuse angle
condition in the \(\mathbb H_K\)-metric. The additional
restriction accounts for the convection contribution. In a
convection-dominated regime its right-hand side can be large on
coarse or highly stretched cells, so effective-tensor alignment alone
need not make the sufficient spatial condition attainable. Substantial
crosswind resolution may be required. No universal crosswind-mesh
threshold follows from this bound without specifying the velocity,
diffusion tensor and element geometry. The cosine-based formulation is dimensionless. The raw minimum of
$\delta_{ij}^{K,H}$ depends on element size and diffusion scale and
should not be used to rank unrelated mesh resolutions. 

The next observation identifies a physical-diffusion restriction implied jointly by the conservative mass and stiffness conditions.

\begin{remark}[A necessary physical-diffusion restriction]
\label{rem:physical-diffusion-restriction}
Suppose the mass restriction \eqref{cond-mij} and the local spatial
condition \eqref{delta-H-cond} both hold for an ordered pair
$i\ne j$ on $K$. Using the definition of $\mathbb H_K$, the latter
condition gives
\[
-\varepsilon(\boldsymbol q_i^K)^T\mathbb D_K\boldsymbol q_j^K
\ge \frac{B_K}{d+1}|\boldsymbol q_j^K|
+\tau_K(\boldsymbol q_i^K)^T
\overline{\boldsymbol b\boldsymbol b^T}_K\boldsymbol q_j^K.
\]
Since the velocity bound implies
$|(\boldsymbol q_i^K)^T
\overline{\boldsymbol b\boldsymbol b^T}_K\boldsymbol q_j^K|
\le B_K^2|\boldsymbol q_i^K||\boldsymbol q_j^K|$,
the mass restriction yields
\begin{equation}\label{eq:necessary-physical-diffusion}
-\varepsilon(\boldsymbol q_i^K)^T\mathbb D_K\boldsymbol q_j^K
\ge
\frac{B_K|\boldsymbol q_j^K|}{(d+1)(d+2)}.
\end{equation}
For $B_K>0$ it follows that the corresponding physical-diffusion
metric angle is strictly acute. Moreover,
$|\boldsymbol q_i^K|=1/h_i^K$ and the spectral bound on
$\mathbb D_K$ give the necessary altitude-based restriction
\begin{equation}\label{eq:necessary-altitude-ratio}
\frac{B_Kh_i^K}{\varepsilon\lambda_{\max}(\mathbb D_K)}
\le (d+1)(d+2).
\end{equation}
In two dimensions the right-hand side is $12$. This restriction is
necessary for the combined conservative elementwise conditions,
not for nodal DMP satisfaction itself and not for the sharper
assembled-matrix criterion. It is not a universal threshold for
unspecified diameter-based mesh P\'eclet numbers or arbitrary
anisotropy. It explains why those sufficient conditions may be
particularly restrictive in convection-dominated regions despite the
presence of streamline diffusion.
\end{remark}

The absolute-value convection estimate in Lemma~\ref{lemma-aij} is
conservative. For a given element, the exact off-diagonal sign can
instead be checked from its signed convection integral. 
The spatial signs must be combined with a time-step restriction to obtain a monotone fully discrete update.

\subsection{Time-step conditions and the \texorpdfstring{$M$}{M}-matrix property}
\label{subsec:temporal-dmp}

For \(i\ne j\), set
\begin{align}
\eta_{ij}^K&:=|K|\delta_{ij}^{K,H},\\
\mu_{ij}^K&:=|K|\left[\frac1{(d+1)(d+2)}+\frac{\tau_KB_K}{d+1}|\boldsymbol q_i^K|\right].
\end{align}
Then the preceding elementwise bounds give, for every
adjacent ordered pair $i\ne j$,
\[
0\le M_{ij}\le\sum_{K\in\omega_i\cap\omega_j}\mu_{ij}^K,
\qquad
A_{ij}\le-\sum_{K\in\omega_i\cap\omega_j}\eta_{ij}^K\le0.
\]
Indeed, the upper estimate for $M_{ij}^K$ uses the positive sign of
the absolute-value bound in the proof of Lemma~\ref{lemma-mij},
whereas the stiffness estimate is Lemma~\ref{lemma-aij}.
Consequently,
\[
M_{ij}+\theta\Delta t A_{ij}
\le\sum_{K\in\omega_i\cap\omega_j}
 (\mu_{ij}^K-\theta\Delta t\eta_{ij}^K).
\]
Assume every denominator appearing in the following adjacent-pair maximum is strictly positive. 
A zero or negative denominator is not assigned an infinite lower-bound value. If \(\theta>0\) and
\begin{equation}
\label{temporallowerbound}
\Delta t\ge\frac1\theta
\max_{i\in I,\,j\ne i:\,\omega_i\cap\omega_j\ne\emptyset}
\frac{\sum_{K\in\omega_i\cap\omega_j}\mu_{ij}^K}
{\sum_{K\in\omega_i\cap\omega_j}\eta_{ij}^K},
\end{equation}
then
\begin{equation}
\label{tildeA-Z}
\widetilde A_{ij}=M_{ij}+\theta\Delta t A_{ij}\le0,
\qquad i\ne j.
\end{equation}

The following lemma identifies the $M$-matrix structure of the Dirichlet-modified time-step operator.

\begin{lemma}[Nonsingular $M$-matrix property]
\label{lemma-effective-mmatrix}
Assume \eqref{cond-mij}, \eqref{delta-H-cond}, and
\eqref{temporallowerbound}. Then the interior principal block of
$\widetilde{\boldsymbol A}=\boldsymbol M+\theta\Delta t\boldsymbol A$
is a nonsingular $M$-matrix. After replacing the Dirichlet rows by
identity equations while retaining the interior-to-boundary columns,
the resulting full matrix has a nonnegative inverse.
\end{lemma}

\begin{proof}
For each interior row $i$, the nodal basis satisfies
$\sum_j\varphi_j=1$ and $\sum_j\nabla\varphi_j=0$ on every element.
Consequently, summing the diffusion, convection and streamline
stiffness entries over all trial indices gives
\[
\sum_jA_{ij}=0,
\qquad
\sum_j\widetilde A_{ij}=\sum_jM_{ij}>0.
\]
Here the sums include boundary trial nodes. In combination with
\eqref{tildeA-Z}, this implies that each interior diagonal is positive
and
\begin{equation}\label{eq:row-margin}
\widetilde A_{ii}-\sum_{j\ne i}|\widetilde A_{ij}|
=\sum_j\widetilde A_{ij}>0.
\end{equation}
Partition the nodes into interior $I$ and boundary $B$ and write
$B_{II}=(\widetilde A_{ij})_{i,j\in I}$ and
$B_{IB}=(\widetilde A_{ij})_{i\in I,j\in B}$.
Because the boundary off-diagonal entries are nonpositive,
\[
(B_{II})_{ii}-\sum_{\substack{j\in I\\j\ne i}}
|(B_{II})_{ij}|
\ge\sum_j\widetilde A_{ij}>0.
\]
Thus $B_{II}$ is a strictly row-diagonally-dominant $Z$-matrix
(a matrix whose off-diagonal entries are nonpositive).
Write $B_{II}=D-N$, where $D$ is its positive diagonal matrix and
$N\ge0$. Strict row dominance gives
\[
\|D^{-1}N\|_\infty
=\max_{i\in I}\frac{\sum_{j\in I,j\ne i}N_{ij}}{D_{ii}}<1.
\]
The Neumann series is therefore convergent and termwise nonnegative:
\begin{equation}\label{eq:m-matrix-inverse}
B_{II}^{-1}
=(I-D^{-1}N)^{-1}D^{-1}
=\sum_{k=0}^\infty(D^{-1}N)^kD^{-1}\ge0.
\end{equation}
This proves both nonsingularity and the $M$-matrix property without
assuming symmetry of the SUPG operator.

With the boundary rows replaced by identity equations, the complete
matrix is
\[
\widehat B=\begin{pmatrix}B_{II}&B_{IB}\\0&I\end{pmatrix}.
\]
Since $B_{IB}\le0$, direct block multiplication yields
\[
\widehat B^{-1}=
\begin{pmatrix}
B_{II}^{-1}&-B_{II}^{-1}B_{IB}\\0&I
\end{pmatrix}\ge0.
\]
The interior-to-boundary coupling is essential for comparison with
nonhomogeneous Dirichlet data. It is not removed from the proof.
\end{proof}

To make
\[
\boldsymbol R:=\boldsymbol M-(1-\theta)\Delta t\boldsymbol A
\]
nonnegative, no additional condition is needed for off-diagonal entries. For diagonal entries define
\begin{equation}
\label{Mii-low}
M_{ii}^{\mathrm low}:=
\sum_{K\in\omega_i}|K|
\left[\frac2{(d+1)(d+2)}-\frac{\tau_KB_K}{d+1}|\boldsymbol q_i^K|\right]>0
\end{equation}
and
\begin{equation}
\label{Aii-up}
A_{ii}^{\mathrm up}:=
\sum_{K\in\omega_i}|K|
\left[\frac{B_K}{d+1}|\boldsymbol q_i^K|
+\lambda_{\max}(\mathbb H_K)|\boldsymbol q_i^K|^2\right].
\end{equation}
Indeed, the diagonal SUPG mass entry is bounded below, element by
element, by
\[
M_{ii}^K
\ge|K|\left[
\frac{2}{(d+1)(d+2)}
-\frac{\tau_KB_K|\boldsymbol q_i^K|}{d+1}
\right],
\]
while the diagonal stiffness entry satisfies
\begin{align*}
A_{ii}^K
&=\int_K\varphi_i(\boldsymbol b\cdot\boldsymbol q_i^K)
\,d\boldsymbol x+|K|(\boldsymbol q_i^K)^T\mathbb H_K\boldsymbol q_i^K\\
&\le|K|\left[\frac{B_K|\boldsymbol q_i^K|}{d+1}
+\lambda_{\max}(\mathbb H_K)|\boldsymbol q_i^K|^2\right].
\end{align*}
Thus $M_{ii}\ge M_{ii}^{\mathrm{low}}>0$ and
$A_{ii}\le A_{ii}^{\mathrm{up}}$. The latter upper bound is strictly
positive when $\varepsilon>0$ and $K$ is nondegenerate. Hence,
for $0\le\theta<1$,
\[
R_{ii}=M_{ii}-(1-\theta)\Delta t A_{ii}
\ge M_{ii}^{\mathrm{low}}
-(1-\theta)\Delta t A_{ii}^{\mathrm{up}}.
\]
It follows that the following explicit sufficient restriction controls
all interior diagonal entries of $R$. For a fixed assembled operator, the
sharper test is simply $M_{ii}-(1-\theta)\Delta t A_{ii}\ge0$ for
all relevant rows. The bound below is a conservative geometric estimate
that avoids assembling the exact diagonal. The assembled-matrix criterion
in Lemma~\ref{lemma:assembled-dmp} tests the actual entries.

Then, for \(0\le\theta<1\),
\begin{equation}
\label{temporalupperbound}
\Delta t\le\frac1{1-\theta}\min_i\frac{M_{ii}^{\mathrm low}}{A_{ii}^{\mathrm up}}
\end{equation}
implies \(R_{ii}\ge0\). When \(\theta=1\), this upper bound disappears because \(\boldsymbol R=\boldsymbol M\), 
whereas the lower bound \eqref{temporallowerbound} remains. The minimum in \eqref{temporalupperbound} 
runs only over interior test nodes with \(A_{ii}^{\mathrm{up}}>0\). 
For $0<\theta<1$, the sufficient DMP conditions require a nonempty interval
between \eqref{temporallowerbound} and \eqref{temporalupperbound}.

For $\theta=0$ (explicit Euler), the left matrix is the consistent
SUPG mass matrix $B=M$. Its off-diagonal entries are generally positive, so
neither the lower-bound argument nor the $M$-matrix theorem below applies.
This does not rule out a DMP for a different mass treatment or under
directly verified update-matrix conditions. For example, when $M_{II}$ is invertible
and the source vanishes, define
\[
W_{\mathrm{old}}^{(0)}=M_{II}^{-1}(M-\Delta t A)_{I,:},
\qquad W_{\mathrm{bdry}}^{(0)}=-M_{II}^{-1}M_{IB}.
\]
If both matrices are entrywise nonnegative, the constant-preserving
row-sum identity gives a nodal interval DMP for the corresponding
explicit update with admissible old and boundary values. These
conditions are not implied by the consistent-mass element estimates here.

\subsection{Bound-preserving solution transfer between adapted meshes}
\label{subsec:remapping-dmp}

For an accepted mesh transfer, the property required for propagating
a common nodal interval is that the transfer map
$\mathcal T^m:\mathbb R^{N_m}\to\mathbb R^{N_{m+1}}$ be
range preserving on that interval:
\begin{equation}\label{eq:range-preserving-transfer}
\boldsymbol v\in[m_*,M_*]^{N_m}
\quad\Longrightarrow\quad
\mathcal T^m(\boldsymbol v)\in[m_*,M_*]^{N_{m+1}}.
\end{equation}
Neither linearity nor conservation is required for this implication.
A standard sufficient realization is barycentric interpolation. Let
\(P^m\) denote the interpolation matrix that maps nodal values on
\(\mathcal T_h^m\) to the nodes of
\(\mathcal T_h^{m+1}\):
\[
\widetilde{\boldsymbol u}^{\,m}=P^m\boldsymbol u^m.
\]
For standard \(P_1\) interpolation, each new node is evaluated in the old simplex containing it. 
Its interpolation weights are barycentric coordinates and therefore are nonnegative and sum to one. Consequently,
\begin{equation}
\label{P-stochastic}
P^m\ge0,
\qquad
P^m\boldsymbol 1=\boldsymbol 1.
\end{equation}
Thus nodal bounds are preserved:
\begin{equation}
\min_j u_j^m\le \widetilde u_i^m\le\max_j u_j^m.
\label{remap-bounds}
\end{equation}
This argument assumes ordinary pointwise \(P_1\) interpolation with a valid containing simplex. 
For a nonlinear transfer, condition~\eqref{eq:range-preserving-transfer}
can instead be checked directly. In particular, a feasible mass
correction with all output vertices constrained to remain in the
prescribed interval is range preserving for inputs in that interval.
This conditional statement does not establish that the particular
numerical correction always has a feasible constrained solve or that
it removes violations already present in the old solution.

The mass, stiffness, time-step, and transfer properties now provide the hypotheses for a discrete comparison result.

\subsection{Discrete monotonicity theorem}

We can now state the one-step comparison theorem, including its conditions on source and boundary data.

\begin{theorem}[A discrete comparison principle for one adaptive SUPG step]
\label{thm-DMP}
Let \(0<\theta\le1\). On the accepted mesh, assume exact assembly of the matrices in \eqref{MmatrixFormula}--\eqref{Amatrixint} 
with the coefficients frozen at the stage time. Suppose that
\begin{enumerate}
\item[(i)] \(\mathbb D_K\) is SPD and \(\tau_K\ge0\) on each nondegenerate element;
\item[(ii)] the mass condition \eqref{cond-mij} holds for every local test vertex;
\item[(iii)] the off-diagonal stiffness condition \eqref{delta-H-cond} holds for every ordered pair of distinct local vertices;
\item[(iv)] the strictly positive denominators and the lower time-step bound in \eqref{temporallowerbound} hold for 
all neighboring pairs with an interior test vertex;
\item[(v)] when \(0<\theta<1\), the upper bound \eqref{temporalupperbound} also holds.
\end{enumerate}
Let \(I\) and \(B\) denote the interior and boundary nodes. Define
\[
B=M+\theta\Delta t A,\qquad R=M-(1-\theta)\Delta t A,
\]
where only the interior test rows are used, and let \(F_I\) denote the assembled stage-weighted load term 
in \eqref{fvector} (before multiplication by \(\Delta t\)). The interior and prescribed-boundary equations are
\begin{equation}
\label{dmp-block-step}
\begin{pmatrix}B_{II}&B_{IB}\\0&I\end{pmatrix}
\begin{pmatrix}u_I^{m+1}\\g_B^{m+1}\end{pmatrix}
=
\begin{pmatrix}R_{II}&R_{IB}\\0&0\end{pmatrix}
\begin{pmatrix}\widetilde u_I^m\\\widetilde u_B^m\end{pmatrix}
+\begin{pmatrix}\Delta t F_I\\g_B^{m+1}\end{pmatrix}.
\end{equation}
The full left matrix in \eqref{dmp-block-step} has a nonnegative inverse, while the full right coefficient matrix 
is nonnegative. Consequently the new nodal solution is nondecreasing in each component of the transferred old-time state, 
in each component of \(F_I\), and in the new boundary data. In particular, for any constants \(a\le b\), if
\[
a\le\widetilde u_j^m\le b,\qquad
 a\le g_j^{m+1}\le b,\qquad F_I=0,
\]
then \(a\le u_j^{m+1}\le b\) for all nodes. For a nonzero source, the lower bound propagates when \(F_I\ge0\), 
and the upper bound propagates when \(F_I\le0\), in each case with the corresponding bound on the old-time and Dirichlet data. 
Both bounds are therefore guaranteed by this theorem for a common fixed interval when \(F_I=0\). 
A general signed source requires a separate source-dependent comparison function.
\end{theorem}

\begin{proof}
For brevity, let $\widehat B$ denote the Dirichlet-modified block
matrix in \eqref{dmp-block-step}. The mass lemma gives $M_{ij}\ge0$
for each interior test row, including boundary trial columns. The
stiffness lemma gives $A_{ij}\le0$ for all distinct pairs. The
lower time-step restriction therefore gives $B_{ij}\le0$ for
$i\in I$, $j\ne i$. Lemma~\ref{lemma-effective-mmatrix} shows that
$B_{II}^{-1}\ge0$, $B_{IB}\le0$, and $\widehat B^{-1}\ge0$.

For $0<\theta<1$, the off-diagonal entries of
$R=M-(1-\theta)\Delta t A$ are nonnegative because $M_{ij}\ge0$
and $A_{ij}\le0$. For diagonal interior entries, the estimates
\eqref{Mii-low}--\eqref{temporalupperbound} imply
\[
R_{ii}\ge M_{ii}^{\mathrm{low}}
-(1-\theta)\Delta t A_{ii}^{\mathrm{up}}\ge0.
\]
For backward Euler ($\theta=1$), $R=M\ge0$ and no upper bound is
needed. It follows that $R_{I,:}=[R_{II}\ R_{IB}]\ge0$ in both cases.

Solving the interior block of \eqref{dmp-block-step} shows the
actual dependence on the old state, new Dirichlet values and load:
\begin{equation}\label{eq:interior-positive-update}
u_I^{m+1}
=\underbrace{B_{II}^{-1}R_{I,:}}_{W_{\mathrm{old}}\ge0}
 \widetilde u^{m}
+\underbrace{(-B_{II}^{-1}B_{IB})}_{W_{\mathrm{bdry}}\ge0}
 g_B^{m+1}
+\Delta t\underbrace{B_{II}^{-1}}_{\ge0}F_I.
\end{equation}
This expression proves componentwise monotonicity in all three data
vectors without invoking any unstated sign property of the source.

Finally, the partition of unity yields, row by row,
\[
A\boldsymbol1=0,
\qquad
B\boldsymbol1=M\boldsymbol1=R\boldsymbol1.
\]
Partition the constant vectors into interior and boundary components.
Multiplication by $B_{II}^{-1}$ gives
\begin{equation}\label{eq:stochastic-update-weights}
W_{\mathrm{old}}\boldsymbol1
+W_{\mathrm{bdry}}\boldsymbol1
=B_{II}^{-1}(R\boldsymbol1-B_{IB}\boldsymbol1)
=\boldsymbol1.
\end{equation}
Thus each row of the combined old-state and new-boundary weights is a
convex combination. If both sets of nodal data lie in $[a,b]$ and
$F_I=0$, \eqref{eq:interior-positive-update}--\eqref{eq:stochastic-update-weights} give
$a\le u_i^{m+1}\le b$ for all interior indices. The prescribed
boundary nodes satisfy the same inequalities by assumption.
For $F_I\ge0$, the source contribution in
\eqref{eq:interior-positive-update} is nonnegative, so the lower bound
remains valid, and analogously $F_I\le0$ suffices for the upper
bound. A source of unrestricted sign requires comparison with an
appropriate nonconstant solution or a source-dependent interval.
\end{proof}

The elementwise bounds above can be conservative. The following lemma gives a separate 
condition based directly on assembled matrix entries.

\begin{lemma}[Direct assembled-matrix DMP criterion]
\label{lemma:assembled-dmp}
For a single accepted time step with $0<\theta\le1$, let $B=M+\theta\Delta t A$
and $R=M-(1-\theta)\Delta t A$ be the actually assembled stage-frozen matrices,
before imposing Dirichlet conditions. Suppose that the interior rows satisfy
\begin{align}
&B_{ij}\le0 \quad(i\in I,\ j\ne i),
\label{eq:direct-B-sign}\\
&\sum_{j\in I\cup\partial I}B_{ij}>0 \quad(i\in I),
\label{eq:direct-row-sum}\\
&R_{ij}\ge0 \quad(i\in I,\ j\in I\cup\partial I),
\label{eq:direct-R-sign}
\end{align}
where $\partial I$ denotes the boundary-node index set, and assume $A\boldsymbol1=0$
on interior rows. Then the Dirichlet-modified left matrix in
\eqref{dmp-block-step} has a nonnegative inverse and the right coefficient
matrix is nonnegative. Hence the comparison and constant-interval conclusions
of Theorem~\ref{thm-DMP} hold, with the same conditions on data and transfer.
These conditions concern the assembled operator, not the individual
elementwise inequalities.
\end{lemma}
\begin{proof}
Let $\widehat B$ be the Dirichlet-modified block matrix, with
$B_{II}$ its interior principal block. For $i\in I$, the sign
assumption $B_{ij}\le0$ for $j\ne i$ and positive full row sum give
\[
B_{ii}=\sum_jB_{ij}+\sum_{j\ne i}|B_{ij}|>0,
\qquad
B_{ii}-\sum_{j\ne i}|B_{ij}|=\sum_jB_{ij}>0.
\]
In particular,
$B_{ii}-\sum_{j\in I,j\ne i}|B_{ij}|>0$.
Consequently $B_{II}=D_B-N_B$, where
$D_B=\operatorname{diag}(B_{II})$ is a positive diagonal matrix and
$N_B=D_B-B_{II}$ is entrywise nonnegative. Strict diagonal dominance gives
$\|D_B^{-1}N_B\|_\infty<1$. The Neumann-series argument in
\eqref{eq:m-matrix-inverse}, applied to $D_B$ and $N_B$, proves
$B_{II}^{-1}\ge0$.
The sign condition also gives $B_{IB}\le0$, hence
\[
\widehat B^{-1}
=\begin{pmatrix}
B_{II}^{-1}&-B_{II}^{-1}B_{IB}\\0&I
\end{pmatrix}\ge0.
\]
The notation $R_{I,:}$ denotes the submatrix formed by all interior
test rows $i\in I$ and all trial columns $j\in I\cup B$.
Thus the assembled right-block condition is $R_{I,:}\ge0$ entrywise.
Since $A\boldsymbol1=0$, the complete row identity
$B\boldsymbol1=R\boldsymbol1$ holds on interior rows.
Thus the two nonnegative weight matrices in
\eqref{eq:interior-positive-update} satisfy the row-sum identity
\eqref{eq:stochastic-update-weights}. The explicit update formula
therefore proves monotonicity and interval preservation, with the
same source and Dirichlet qualifications as in
Theorem~\ref{thm-DMP}. This proof requires neither local
$\delta_{ij}^{K,H}\ge0$ nor elementwise mass inequalities: the
assembled signs are verified directly instead.
\end{proof}

The distinction between local sufficient tests and the assembled-matrix criterion 
is important when interpreting the numerical diagnostics.

\begin{remark}[Local conditions versus algebraic certification]
\label{rem:local-global-sign}
The local inequalities in Theorem~\ref{thm-DMP} imply
\eqref{eq:direct-B-sign}--\eqref{eq:direct-R-sign}, but their converse is
not required and generally does not hold: contributions from neighboring
elements can offset one another in the assembled stiffness entries.
Thus an elementwise failure does not by itself rule out the direct test.
For backward Euler ($\theta=1$), $R=M$, so the latter test still requires
nonnegative relevant entries of the assembled SUPG mass matrix. 
\end{remark}

Applying the one-step comparison result after each range-preserving transfer gives a global-in-time nodal bound.

\begin{corollary}[Homogeneous-source nodal DMP under mesh adaptation]
\label{cor:DMP}
Suppose the hypotheses of Theorem~\ref{thm-DMP} (or the direct
assembled-matrix conditions of Lemma~\ref{lemma:assembled-dmp}) hold
at each accepted time step, $f=0$, and every accepted mesh transfer
$\mathcal T^m$ satisfies \eqref{eq:range-preserving-transfer} for a
common interval $[m_*,M_*]$. If the initial nodal state and all
prescribed boundary values lie in this interval, then
\[
m_*\le u_i^m\le M_*,\qquad m=0,1,\ldots,N_t.
\]
\end{corollary}
\begin{proof}
The assertion follows by induction over accepted physical-time
steps. The initial state is admissible. If
$\boldsymbol u^m\in[m_*,M_*]^{N_m}$, range preservation gives
$\widetilde{\boldsymbol u}^m=\mathcal T^m(\boldsymbol u^m)\in[m_*,M_*]^{N_{m+1}}$. 
The homogeneous-source SUPG load vanishes.
The matrix comparison theorem, together with admissible next-time
Dirichlet data, then yields
$\boldsymbol u^{m+1}\in[m_*,M_*]^{N_{m+1}}$.
The predictor--remesh iterations do not advance physical time and
do not enter the induction as accepted states.
\end{proof}

The theorem supplies sufficient, not necessary, conditions for
bound preservation. In particular, failure of a local condition does
not imply an actual violation. The sign of a stabilized source vector
must also be checked: a nonnegative continuous source does not by
itself make the SUPG load vector nonnegative.

We use the effective diffusion tensor identified above to design metrics that 
favor the anisotropic geometry appearing in the DMP conditions.

\subsection{DMP-oriented and adaptive metrics for anisotropic remeshing}
\label{subsec-metric}

The tensor entering the off-diagonal matrix condition determines a natural geometric reference metric,
\begin{equation}
\label{Meff}
\mathbb{M}_{H,K}=c_{H,K}\mathbb H_K^{-1},\qquad c_{H,K}>0.
\end{equation}
Since \(\mathbb H_K\) is SPD, \eqref{Meff} is well defined without an additional definiteness assumption 
involving derivatives of the velocity. The inverse tensor fixes the principal directions and anisotropy ratio of the mesh, 
whereas the positive scalar \(c_{H,K}\) remains available to control local mesh density. 
A convenient choice, in two dimensions, is
\begin{equation}\label{eq:normalized-pure-H}
c_{H,K}=\det(\mathbb H_K)^{1/2},\Rightarrow
\mathbb{M}_{H,K}=\det(\mathbb H_K)^{1/2}\mathbb H_K^{-1},
\quad \text{ and } \det(\mathbb{M}_{H,K})=1.
\end{equation}
This element-stage determinant normalization supplies constant input metric
density while retaining the effective-tensor eigenvectors and anisotropy
ratios. Subsequent area-weighted interpolation to vertices, global complexity
scaling, and remeshing can alter the local determinant. It does not prescribe
uniform physical triangle areas. A mesh approximately uniform in this metric tends to make elements closer 
to isotropic in the \(\mathbb H_K\)-metric and therefore promotes the metric-angle part of \eqref{effective-angle-condition}. 
Nevertheless, the metric condition alone cannot guarantee a DMP. The convection-inclusive restriction \eqref{delta-H-cond} 
and the physical-diffusion limitation in Remark~\ref{rem:physical-diffusion-restriction} must still be considered.

The free scalar in \eqref{Meff} can be used for interpolation-error adaptation without changing the \(\mathbb H_K^{-1}\) alignment. 
This DMP-oriented construction differs from the anisotropic SUPG metric of Di, Xie, and Yin \cite{di2018anisotropic}, 
whose design uses the solution Hessian and a stabilized discretization-error estimate rather than 
the sign conditions of the present fully discrete operator. In particular, their notation for the solution Hessian 
should not be confused with our effective diffusion tensor \(\mathbb H_K\). We follow the construction of 
Li and Huang \cite{li2010anisotropic}, who considered the family \(\theta_K\mathbb D_K^{-1}\) for 
anisotropic diffusion and selected \(\theta_K\) by minimizing a regularized \(H^1\)-interpolation-error bound. 
The following estimate applies to a frozen SPD tensor under ideal metric alignment. It motivates a density prescription 
but does not establish optimality of the coupled SUPG computation. The spectral and Frobenius matrix norms are 
equivalent in fixed dimension. Consequently using the spectral norm modifies the bound by a dimension-dependent constant 
rather than changing its asymptotic order. Let \(H(u)=\nabla^2u\) and let \(|H(u)|\) denote the matrix absolute value 
obtained by replacing the Hessian eigenvalues by their absolute values. Define
\begin{equation}
\label{Psi-A}
\Psi_{A,K}
:=
\det(\mathbb A_K)^{-1/d}
\|\mathbb A_K^{-1}\|_2
\frac{1}{|K|}
\int_K
\|\mathbb A_K |H(u)|\|_2^2\,d\boldsymbol x.
\end{equation}
Then we have the following lemma.

\begin{lemma}[Conditional anisotropic interpolation estimate]
\label{lemma:general-spd-interpolation}
Let $u\in H^2(\Omega)$ in two space dimensions and let
$\mathbb A_K$ be a constant SPD matrix on each nondegenerate simplex
$K$. Let $G_K=F'_K$ be the Jacobian of the affine map from a fixed
unit-volume regular reference simplex $\widehat K$ to $K$.
Assume the ideal alignment relation
\begin{equation}\label{eq:ideal-A-alignment}
G_K^T\mathbb A_K^{-1}G_K=\sigma_K I,
\qquad\sigma_K>0.
\end{equation}
Then, with $\Psi_{A,K}$ defined in \eqref{Psi-A},
\begin{equation}
\label{A-interp-bound}
|u-\Pi_hu|_{H^1(\Omega)}^2
\le C\sum_{K\in\mathcal T_h}|K|^{(d+2)/d}\Psi_{A,K},
\end{equation}
where $C$ depends on the fixed reference simplex and dimension,
not on the element aspect ratios or the eigenvalues of
$\mathbb A_K$.
\end{lemma}

\begin{proof}
Let $\widehat u=u\circ F_K$ and denote by
$\widehat\Pi$ the $P_1$ nodal interpolant on $\widehat K$.
The standard reference-simplex interpolation inequality
and the affine change of coordinates give
\begin{align*}
\|\nabla(u-\Pi_hu)\|_{L^2(K)}^2
&\le\|G_K^{-1}\|_2^2|\det G_K|
\|\nabla_{\widehat x}(\widehat u-\widehat\Pi\widehat u)
\|_{L^2(\widehat K)}^2\\
&\le C\|G_K^{-1}\|_2^2
\int_K\|G_K^T H(u)G_K\|_F^2\,d\boldsymbol x.
\end{align*}
Since \eqref{eq:ideal-A-alignment} holds, the matrix
$O_K=\sigma_K^{-1/2}\mathbb A_K^{-1/2}G_K$ is orthogonal,
and hence
\[
G_K=\sigma_K^{1/2}\mathbb A_K^{1/2}O_K,
\qquad
\|G_K^{-1}\|_2^2=\sigma_K^{-1}\|\mathbb A_K^{-1}\|_2.
\]
Moreover, $H(u)$ is symmetric almost everywhere. The symmetric
matrix $\mathbb A_K^{1/2}H(u)\mathbb A_K^{1/2}$ is similar to
$\mathbb A_K H(u)$, so its spectral norm equals its spectral radius,
which is no larger than $\|\mathbb A_K H(u)\|_2$. The Frobenius norm
is at most $\sqrt d$ times the spectral norm. Thus
\[
\|G_K^TH(u)G_K\|_F^2
\le d\sigma_K^2\|\mathbb A_KH(u)\|_2^2.
\]
Since $H(u)=|H(u)|S(u)$ for an orthogonal symmetric sign matrix
$S(u)$ obtained from the Hessian's spectral decomposition,
$\|\mathbb A_KH(u)\|_2=\|\mathbb A_K|H(u)|\|_2$.
Taking determinants in \eqref{eq:ideal-A-alignment}, using
$|\det G_K|=|K|$ for the unit-volume reference simplex, gives
\[
\sigma_K
=|K|^{2/d}\det(\mathbb A_K)^{-1/d}.
\]
Substitution into the reference estimate now yields
\begin{align*}
\|\nabla(u-\Pi_hu)\|_{L^2(K)}^2
&\le C |K|^{2/d}\det(\mathbb A_K)^{-1/d}
\|\mathbb A_K^{-1}\|_2
\int_K\|\mathbb A_K|H(u)|\|_2^2\,d\boldsymbol x\\
&=C|K|^{(d+2)/d}\Psi_{A,K}.
\end{align*}
Summing over all elements proves \eqref{A-interp-bound}.
\end{proof}

Lemma~\ref{lemma:general-spd-interpolation} is conditional on exact alignment. The subsequent metric
formula uses the resulting estimate as a design criterion.
The numerical meshes need not satisfy \eqref{eq:ideal-A-alignment}
exactly, and the estimate is not a convergence theorem for the
SUPG-discretized PDE. In particular, $\mathbb A_K^{-1}$ fixes the
principal directions and eigenvalue ratios of the mesh metric. The scalar
interpolation factor changes element density but cannot rotate those
directions to match an unrelated solution layer.
Regularizing with a positive parameter \(\alpha_{A,h}\) gives
\begin{equation}
\label{A-interp-reg}
|u-\Pi_hu|_{H^1(\Omega)}^2
\le C\alpha_{A,h}
\sum_{K\in\mathcal T_h}
|K|^{(d+2)/d}
\left(1+\frac{\Psi_{A,K}}{\alpha_{A,h}}\right).
\end{equation}
To explain the equidistribution rule, put
$p=d/(d+2)\in(0,1)$ and
$a_K=|K|(\alpha_{A,h}+\Psi_{A,K})^p\ge0$.
Convexity of $x\mapsto x^{1/p}$ (equivalently, the finite-dimensional
H\"older inequality) gives
\begin{align*}
\sum_K |K|^{1/p}(\alpha_{A,h}+\Psi_{A,K})
&=\sum_K a_K^{1/p}\\
&\ge N^{1-1/p}\left(\sum_K a_K\right)^{1/p}\\
&=N^{-2/d}\left[
\sum_K |K|(\alpha_{A,h}+\Psi_{A,K})^p
\right]^{1/p}.
\end{align*}
For a fixed element count $N$ and a fixed total
$S=\sum_K a_K$, strict convexity of $x^{1/p}$ means that
$\sum_K a_K^{1/p}$ is minimized subject to
$\sum_K a_K=S$ precisely when $a_K=S/N$ for every element.
This is a finite-dimensional constrained statement, not an
unconstrained minimization over meshes: $\Psi_{A,K}$, element areas,
$S$, and even $N$ can change during remeshing. Accordingly, its
equality case identifies an equidistribution target under
fixed total indicator mass and fixed complexity, rather than an
unconditional minimizer of the SUPG discretization error. The subsequent metric is a design rule
under ideal alignment and a fixed-complexity normalization. Comparing this condition with 
metric equidistribution identifies the metric volume density (up to its overall normalization) 
as \(q_{A,K}=(1+\Psi_{A,K}/\alpha_{A,h})^p\). Thus, under the assumed metric alignment, the equality case in this
lower-bound/equidistribution design argument motivates the choice
\begin{equation}
\label{q-A}
q_{A,K}:=\sqrt{\det(\mathbb{M}_{A,K})}
=\left(1+\frac{\Psi_{A,K}}{\alpha_{A,h}}\right)^{d/(d+2)}.
\end{equation}
If \(\mathbb{M}_{A,K}=\theta_{A,K}\mathbb A_K^{-1}\), then
\[
q_{A,K}=\theta_{A,K}^{d/2}\det(\mathbb A_K)^{-1/2}.
\]
Combining this identity with \eqref{q-A} gives
\begin{equation}
\label{M-A-adapt}
\mathbb{M}_{A+\mathrm{adapt},K}
=
\left(1+\frac{\Psi_{A,K}}{\alpha_{A,h}}\right)^{2/(d+2)}
\det(\mathbb A_K)^{1/d}\mathbb A_K^{-1}.
\end{equation}
Following \cite{li2010anisotropic}, we choose the regularization parameter as
\begin{equation}
\label{alpha-A}
\alpha_{A,h}
=
\left[
\frac{1}{|\Omega|}
\sum_{K\in\mathcal T_h}|K|\,
\Psi_{A,K}^{d/(d+2)}
\right]^{(d+2)/d}.
\end{equation}
For nonzero indicators, the normalization has a simple bound,
but should not be interpreted as an equality constraint. With
$p=d/(d+2)$, concavity and subadditivity give
$(1+x)^p\le1+x^p$ for $x\ge0$. Therefore
\begin{align*}
\sum_K |K|q_{A,K}
&\le |\Omega|+\alpha_{A,h}^{-p}
\sum_K|K|\Psi_{A,K}^p\\
&=|\Omega|+|\Omega|=2|\Omega|,
\end{align*}
where the last equality is the definition \eqref{alpha-A}.
Thus this particular choice of $\alpha_{A,h}$ gives an upper
bound on the metric volume, not generally equality. An independent
global metric-complexity normalization can still be applied.

If \(\Psi_{A,K}=0\) for every element, \eqref{alpha-A} gives zero and the regularized quotient is undefined. 
We therefore define the limiting metric \(\det(\mathbb A_K)^{1/d}\mathbb A_K^{-1}\) directly in that degenerate case 
(and normalize the global complexity separately).

Applying \eqref{M-A-adapt} with \(\mathbb A_K=\mathbb H_K\) defines the SUPG-aware DMP-oriented adaptive metric
\begin{equation}
\label{M-H-adapt}
\mathbb{M}_{H+\mathrm{adapt},K}
=
\left(1+\frac{\Psi_{H,K}}{\alpha_{H,h}}\right)^{2/(d+2)}
\det(\mathbb H_K)^{1/d}\mathbb H_K^{-1},
\end{equation}
where
\begin{equation}
\label{Psi-H}
\Psi_{H,K}
=
\det(\mathbb H_K)^{-1/d}
\|\mathbb H_K^{-1}\|_2
\frac{1}{|K|}
\int_K
\|\mathbb H_K|H(u)|\|_2^2\,d\boldsymbol x,
\end{equation}
and \(\alpha_{H,h}\) is given by \eqref{alpha-A} with \(\Psi_{A,K}=\Psi_{H,K}\). In two dimensions this simplifies to
\begin{equation}
\label{M-H-adapt-2d}
\mathbb{M}_{H+\mathrm{adapt},K}
=
\left(1+\frac{\Psi_{H,K}}{\alpha_{H,h}}\right)^{1/2}
\det(\mathbb H_K)^{1/2}\mathbb H_K^{-1},
\end{equation}
with
\[
\alpha_{H,h}
=
\left[
\frac{1}{|\Omega|}
\sum_{K\in\mathcal T_h}|K|\Psi_{H,K}^{1/2}
\right]^2.
\]
In computation, the exact Hessian in \eqref{Psi-H} is replaced by the recovered Hessian described in Section~\ref{subsec:hessian-adapt}. 
Within each BAMG adaptation iteration, \(\mathbb H_K\) is evaluated on the current candidate mesh, then \(\Psi_{H,K}\) 
and the metric are formed. After remeshing, \(\tau_K\), \(\mathbb H_K\), and the metric are recomputed if another adaptation iteration 
is required. Thus the direct metric is always associated with the geometry on which it is evaluated.

For comparison, setting \(\mathbb A_K=\mathbb D_K\) in \eqref{M-A-adapt} recovers the Li--Huang diffusion-based metric
\begin{equation}
\label{MDMP-adapt}
\mathbb{M}_{\mathrm{DMP+adapt},K}
=
\left(1+\frac{\Psi_{D,K}}{\alpha_{D,h}}\right)^{2/(d+2)}
\det(\mathbb D_K)^{1/d}\mathbb D_K^{-1},
\end{equation}
with \(\Psi_{D,K}\) and \(\alpha_{D,h}\) defined by \eqref{Psi-A} and \eqref{alpha-A}. 
The diffusion-only DMP metric without interpolation adaptation is the family
\begin{equation}
\label{MDMP}
\mathbb{M}_{\mathrm{DMP},K}=c_{D,K}\mathbb D_K^{-1},\qquad c_{D,K}>0.
\end{equation}
The derived metric \eqref{M-H-adapt} therefore differs from \eqref{MDMP-adapt} by replacing the physical diffusion geometry with 
the SUPG-effective geometry in both the alignment tensor and the interpolation indicator. More precisely, 
$\mathbb{M}_{H+\mathrm{adapt}}$ and $\mathbb{M}_{\mathrm{DMP+adapt}}$ use the same functional form, $\Psi$ construction, 
and regularization rule, substituting $\mathbb H_K$ for $\mathbb D_K$. Their comparison thus tests this tensor substitution 
within a common construction, rather than a separate normalization formula. It does not isolate tensor orientation alone, 
since the substitution also changes the interpolation indicator and the resulting scalar density. 
In contrast, the adaptation-only metric $\mathbb{M}_{\mathrm{adapt}}$ uses a different construction. In particular,
$\mathbb{M}_{\mathrm{adapt}}$ uses the regularized recovered-Hessian
matrix $\mathbb A_K=I+|H_K(u_h)|/\alpha_h$ and the spectral-norm
density $\rho_K$ in \eqref{adapt-rhoK}, normalized by
\eqref{alpha-adapt-normalization}. The diffusion- and effective-tensor
metrics use the inverse-$\mathbb A_K$ alignment family and the
$\Psi_{A,K}$ moment normalization \eqref{alpha-A}.
Consequently, there is no claim that \eqref{metrictensordef} is
obtained by substituting the regularized Hessian into
\eqref{M-A-adapt}. These are different interpolation-based metric
constructions, not algebraically identical formulas.

The construction is consistent with the diffusion-based Li--Huang
metric. If $\tau_K=0$, then $\mathbb H_K=\varepsilon\mathbb D_K$.
The definition~\eqref{Psi-A} is unchanged by replacing
$\mathbb A_K$ by a positive scalar multiple, so
$\Psi_{H,K}=\Psi_{D,K}$ and the corresponding regularization
parameters agree. Since
$\det(\varepsilon\mathbb D_K)^{1/d}
(\varepsilon\mathbb D_K)^{-1}=
\det(\mathbb D_K)^{1/d}\mathbb D_K^{-1}$,
we obtain $\mathbb{M}_{H+\mathrm{adapt},K}=
\mathbb{M}_{\mathrm{DMP+adapt},K}$ in this limit.

For pure diffusion, alignment with $\mathbb D_K^{-1}$ is directly
related to the anisotropic nonobtuse-angle condition
\cite{li2010anisotropic}. For convection--diffusion, effective-tensor
alignment is only one part of the sufficient condition. We therefore
refer to the resulting metric as DMP-oriented rather than
unconditionally DMP-preserving.

\section{Numerical results}
\label{sec:numerical-experiments}
We compare fixed-mesh SUPG with $\mathbb{M}_{\mathrm{adapt}}$,
$\mathbb{M}_{\mathrm{DMP+adapt}}$, determinant-normalized $\mathbb{M}_H$,
and $\mathbb{M}_{H+\mathrm{adapt}}$. The principal transport tests use
mass-limited stabilization, common time steps and adaptation schedules,
and sign-aware transfer. $N$ denotes the nominal number of triangles.
BAMG does not preserve the element count exactly. The manufactured
solutions use nominal SUPG, whereas transport errors are evaluated
against independent numerical references. Metric comparisons are made
at matched nominal budgets, not identical meshes.

For a target interval $[m_*,M_*]$ and an accepted physical-time
state $u_h(t_n)$, define the nodal bound violations
$E_-(t_n)=\max_i(m_*-u_i^n)_+$ and
$E_+(t_n)=\max_i(u_i^n-M_*)_+$. We report
$E_\pm^{\max}=\max_{0\le n\le N_t}E_\pm(t_n)$, together with the
following whole-trajectory lower-bound measures for the examples with
$m_*=0$:
\begin{equation}\label{eq:whole-time-violations}
\begin{aligned}
E_-^{\max}&=\max_{0\le n\le N_t}E_-(t_n),\\
I_-&=\sum_{n=0}^{N_t-1}\frac{\Delta t}{2}\bigl(E_-(t_n)+E_-(t_{n+1})\bigr),\\
I_{-,2}&=\sum_{n=0}^{N_t-1}\frac{\Delta t}{2}\bigl(\|u_h^-(t_n)\|_{L^2(\Omega)}+\|u_h^-(t_{n+1})\|_{L^2(\Omega)}\bigr),
\end{aligned}
\end{equation}
The integrals are composite trapezoidal approximations in physical time.
$\|u_h^-\|_{L^2}$ is integrated exactly over the $P_1$ sign polygons.
Only accepted post-PDE states are included: these statistics do not bound
violations between recorded levels. Corollary~\ref{cor:DMP} guarantees bound preservation only when its
hypotheses hold at every accepted time step. As noted in
Remark~\ref{rem:local-global-sign}, failure of a sufficient matrix-sign
test does not imply an actual bound violation.
Pure $\mathbb{M}_H$ is retained as a geometric reference that
emphasizes effective-tensor alignment without the recovered-Hessian
density factor. This comparison does not isolate all effects of scalar
metric scaling or establish DMP or accuracy optimality.

We first examine the approximation errors of the five metric choices on a smooth manufactured solution.

\subsection{Example 1: manufactured-solution accuracy and convergence}
\label{subsec:accuracy-benchmark}
On \(\Omega=(0,1)^2\) we use
\begin{align}
 u(x,y,t)&=16\sin(\pi t)x(1-x)y(1-y)\notag\\
 &\quad\times\left[\frac12+\frac1\pi\arctan\!\left(2\varepsilon^{-1/2}
 \left(0.25^2-(x-0.5)^2-(y-0.5)^2\right)\right)\right].
\end{align}
For this accuracy test we set
\[
\mathbb D(\boldsymbol x)=\mathbb I,\qquad \varepsilon=0.1,\qquad \boldsymbol b=(2,3)^T.
\]
The source and boundary data follow from the exact solution:
\begin{equation}
\label{ex1-manufactured-data}
f=u_t-0.1(u_{xx}+u_{yy})+2u_x+3u_y,\qquad
u_0(\boldsymbol x)=0,\qquad g=u|_{\partial\Omega}=0,
\end{equation}
with all derivatives evaluated from the displayed exact solution. We take $T=0.1$, $\Delta t=0.005$, and Crank--Nicolson time
integration, with remeshing every $0.025$ physical-time units (every five PDE steps). 
This fixed schedule reduces remeshing and transfer overhead and was not optimized; 
remeshing after every step is not required by the method.
Convergence is plotted against the actual final triangle count.
This example assesses accuracy on a smooth solution, rather than convection-dominated bound violations.

\begin{figure}[htbp]
\centering
\subfloat[$L^2$ error]{\manuscriptfigure{.48\textwidth}{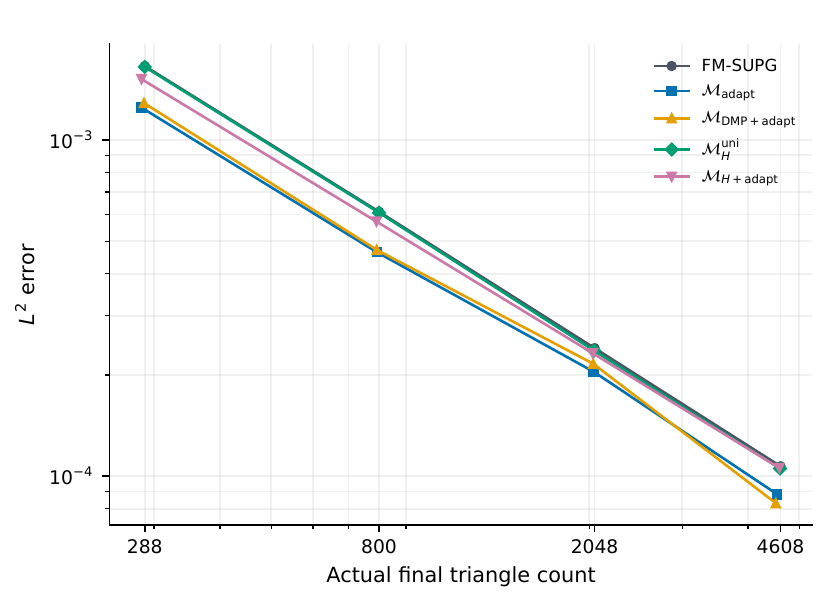}}\hfill
\subfloat[$H^1$ seminorm error]{\manuscriptfigure{.48\textwidth}{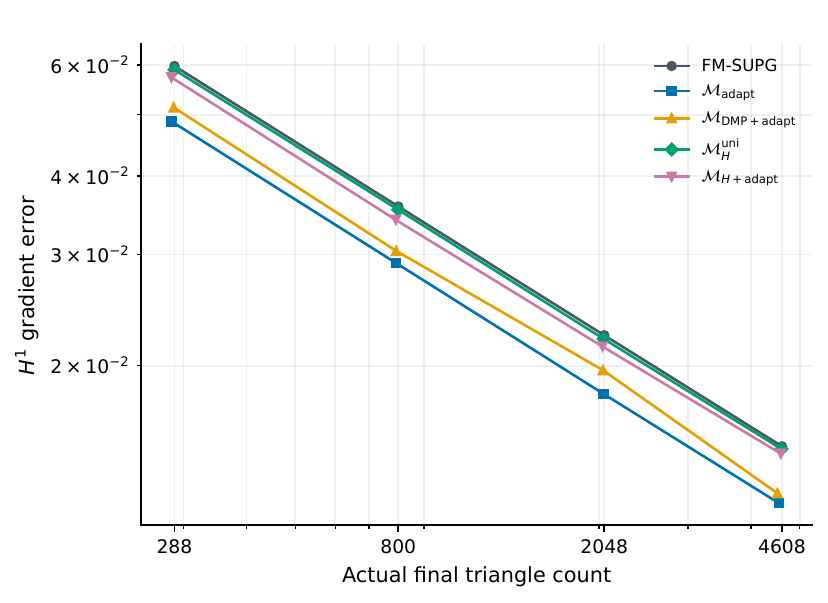}}
\caption{Example~1: $L^2$ and $H^1$-seminorm errors for the five methods versus actual final triangle counts, 
using Crank--Nicolson time integration with $\Delta t=0.005$. The lines $N^{-1}$ and $N^{-1/2}$ are reference slope guides.}
\label{fig:ex1-spatial}
\end{figure}

Figure~\ref{fig:ex1-spatial} shows decreasing errors broadly consistent with the expected $P_1$ rate guides. 
Hessian-only adaptation has the smallest measured $H^1$ error at every nominal target and the smallest $L^2$ error 
at targets $288$, $800$ and $2048$. At target $4608$, $\mathbb{M}_{\mathrm{DMP+adapt}}$ has 
the smallest $L^2$ error ($8.3113\times10^{-5}$ versus $8.8281\times10^{-5}$ for $\mathbb{M}_{\mathrm{adapt}}$), 
but its $H^1$ error is larger. These comparisons concern nominally matched budgets with slightly different actual final counts. 
The distinct goals of interpolation-error adaptation and DMP-oriented geometry do not imply a common accuracy ranking.

We next assess practical bound violations in a convection-dominated problem with isotropic physical diffusion.

\subsection{Example 2: isotropic incompressible Gaussian-hill transport}
\label{subsec:isotropic-dmp}
This convection-dominated test examines numerical bound violations and
the accompanying accuracy cost for isotropic physical diffusion.
Let $\Omega=(0,10)\times(-1.25,1.25)$,
$\varepsilon=10^{-3}$, $\mathbb D=I$, $f=0$, and impose homogeneous
Dirichlet boundary condition. The prescribed initial profile is
\[
 u_0=e^{-((x-1)^2+y^2)/0.25}
 +e^{-((x-1)^2+(y-0.5)^2)/0.25}
 +e^{-((x-1)^2+(y+0.5)^2)/0.25},
\]
with upper target bound $M_0=1+2/e$. The divergence-free transport field is
\[
\boldsymbol b=\left(1-e^{\beta x}\cos(2\pi y),
\frac{\beta}{2\pi}e^{\beta x}\sin(2\pi y)\right),\quad
\beta=\frac1{2\varepsilon}-\sqrt{\left(\frac1{2\varepsilon}\right)^2+4\pi^2}.
\]
We use $T=0.5$, backward Euler with $\Delta t=0.025$, two target
complexities (1152 and 4608 triangles), and mass-compatible limited SUPG.
The analytic Gaussian profile has a nonzero trace on parts of
$\partial\Omega$, whereas the prescribed Dirichlet condition is zero for
$t>0$. We solve this initial-boundary value problem using the Gaussian
profile at interior initial nodes and impose zero values at boundary
nodes. The modification concerns the boundary trace of the initial
finite element field, not an alternative interior initial condition.
The continuous initial and boundary traces are not classically compatible
at $t=0$, but the homogeneous-Dirichlet parabolic problem remains
well defined in a weak sense. This produces an early transient that may
affect reference comparisons. Backward-Euler damping also influences
the computed undershoot amplitudes.

Accuracy is evaluated against an independently refined
Rannacher--Crank--Nicolson reference with $3{,}838{,}464$ triangles
and $\Delta t=0.0015625$. The last spatial refinement produced
relative changes of $0.0112\%$ in $L^2$ and $1.44\%$ in the
$H^1$ seminorm, so the latter reference comparisons are indicative.
The reference's sampled undershoot is $8.66\times10^{-4}$ at
$t\approx0.00469$, whereas the adaptive maxima occur between
$t=0.3$ and $0.5$. The comparison is specific
to the stated final time and initial-boundary data. The results are
presented in Table~\ref{tab:e2-v3-local}.

\begin{table}[htbp]
\centering\scriptsize\setlength{\tabcolsep}{2.5pt}
\begin{tabular}{lrrrrrrrr}
\toprule
Metric & $N_T$ & $E_-^{\max}$ & $I_-$ & $I_{-,2}$ & $E_+^{\max}$ & $L^2_{\mathrm{ref}}(\%)$ & $H^1_{\mathrm{ref}}(\%)$ & Peak/ref $(\%)$\\
\midrule
\multicolumn{9}{l}{Target 1152 triangles}\\
FM--SUPG & 1152 & 0.089741 & 0.013283 & 0.0030447 & 0.15973 & 16.002 & 60.481 & 109.71\\
$\mathbb{M}_{\mathrm{adapt}}$ & 1146 & 0.018812 & 0.0027738 & 0.00050165 & 0 & 12.853 & 31.357 & 97.157\\
$\mathbb{M}_{\mathrm{DMP+adapt}}$ & 1160 & 0.034451 & 0.0037525 & 0.00070537 & 0 & 14.198 & 38.692 & 96.545\\
${\mathbb{M}_H}$ & 1152 & 0.116872 & 0.02194 & 0.004151 & 0 & 22.97 & 57.675 & 92.795\\
$\mathbb{M}_{H+\mathrm{adapt}}$ & 1156 & 0.00291 & 0.00062516 & 0.00011493 & 0 & 20.238 & 43.505 & 90.192\\
\multicolumn{9}{l}{Target 4608 triangles}\\
FM--SUPG & 4608 & 0.0069033 & 0.0011254 & 0.00011771 & 0.021851 & 8.5169 & 35.71 & 101.73\\
$\mathbb{M}_{\mathrm{adapt}}$ & 4610 & 0.0027862 & 0.00033835 & 2.8195e-05 & 0 & 8.8392 & 19.433 & 98.934\\
$\mathbb{M}_{\mathrm{DMP+adapt}}$ & 4638 & 0.0041463 & 0.00051714 & 4.9447e-05 & 0 & 9.4182 & 23.357 & 99.129\\
${\mathbb{M}_H}$ & 4623 & 0.00448618 & 0.0010655 & 0.00011634 & 0 & 9.8414 & 28.781 & 98.571\\
$\mathbb{M}_{H+\mathrm{adapt}}$ & 4613 & 0.0011553 & 0.00010458 & 6.9869e-06 & 0 & 10.102 & 22.46 & 98.688\\
\bottomrule\end{tabular}
\caption{Example~2 with $\Delta t=0.025$ and sign-aware, globally mass-preserving transfer. 
All five methods use the paired principal protocol. Pure $\mathbb{M}_H$ is determinant-normalized in two dimensions. 
Whole-time bound measures use accepted levels and exact spatial sign integration. 
Errors are relative to a qualified numerical reference, not an exact solution. $H^1$ is the gradient seminorm.}
\label{tab:e2-v3-local}
\end{table}

The combined metric $\mathbb{M}_{H+\mathrm{adapt}}$ has the smallest $E_-^{\max}$ at both budgets.
Compared with pure normalized $\mathbb{M}_H$, its reductions are
$97.51\%$ and $74.25\%$. Relative to $\mathbb{M}_{\mathrm{adapt}}$
they are $84.5\%$ and $58.5\%$. Pure $\mathbb{M}_H$ is worse
than fixed SUPG at the coarse budget but better at the fine budget.
The combined metric incurs greater reference error: at target $4608$
its relative $L^2$ error is $10.10\%$, compared with $8.84\%$ for
Hessian-only adaptation. The reported upper-bound excess is nonzero
for fixed SUPG and zero for all adaptive methods. These observed
violation reductions do not constitute DMP certification.

\begin{figure}[!htbp]
\centering
\subfloat[$\mathbb{M}_{\mathrm{adapt}}$, $N_T=1146$]{%
  \manuscriptfigure{.48\textwidth}{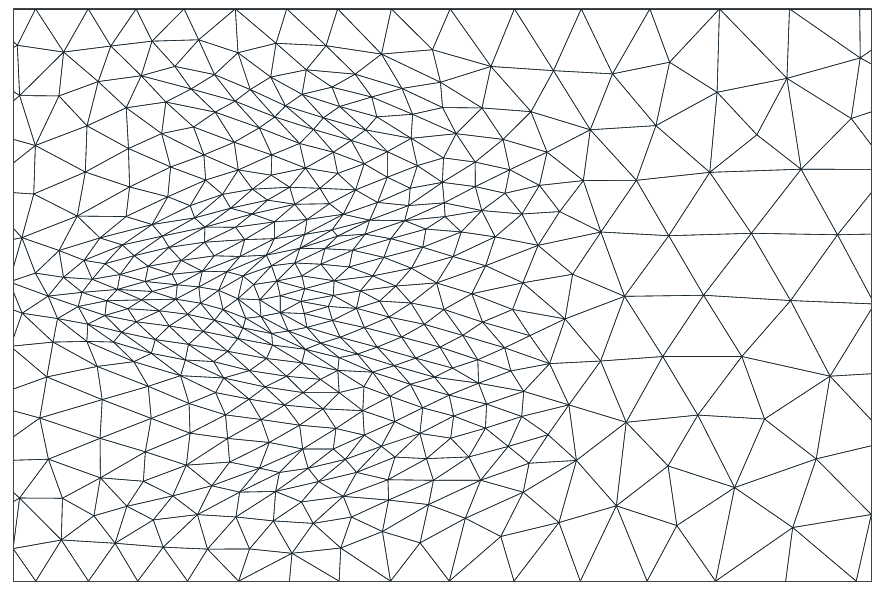}}\hfill
\subfloat[$\mathbb{M}_{H+\mathrm{adapt}}$, $N_T=1156$]{%
  \manuscriptfigure{.48\textwidth}{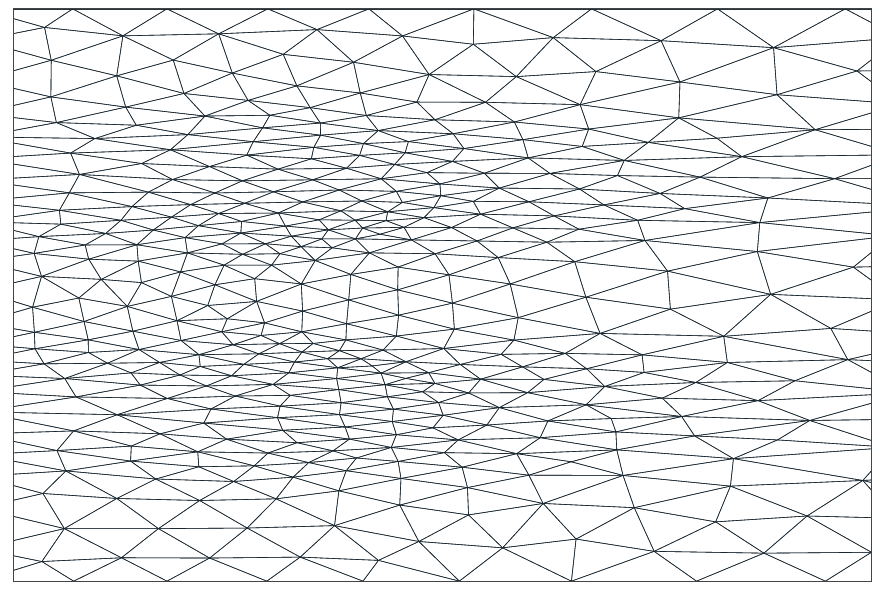}}
\caption{Example~2: final adapted triangular meshes at $T=0.5$ and
nominal target $N\approx1152$. Both panels show the same enlarged
region $0.15\le x\le3.9$, $-1.25\le y\le1.25$ with equal coordinate
scaling. The displayed meshes correspond to the accepted solutions
reported in Table~\ref{tab:e2-v3-local}.}
\label{fig:e2-mesh-pair}
\end{figure}

Figure~\ref{fig:e2-mesh-pair} compares the final adapted meshes at the
coarser target. The two metrics produce different element orientations
and distributions near the transported Gaussian hills.

The next manufactured example examines accuracy when the diffusion tensor is strongly anisotropic.

\subsection{Example 3: effect of solution-dependent metric density}
\label{subsec:anisotropic-accuracy}
To examine the role of solution-dependent density in an anisotropic setting, we use
on $\Omega=(0,1)^2$ the manufactured moving layer
\[
u(x,y,t)=\frac{1}{1+\exp(C(x+y-t))},\qquad C=6,
\]
with $\varepsilon=10^{-3}$,
\[
\mathbb D=\begin{pmatrix}50.5&49.5\\49.5&50.5\end{pmatrix},\qquad
\boldsymbol b=(x+t,y+t).
\]
The eigenvalues of $\mathbb D$ are $100$ and $1$; the source and
Dirichlet data are determined by the exact solution. We use
$T=0.3$, Crank--Nicolson time integration with $\Delta t=0.005$,
and remesh every $0.05$ physical-time units (every ten PDE steps).
This fixed schedule limits repeated mesh generation and transfer and
was not optimized.

We evaluated $L^2$ and $H^1$-seminorm errors for all five methods
at nominal targets $N=288$, $800$, $2048$, and $4608$. Both errors
decrease with refinement for each method; $\mathbb M_{\mathrm{adapt}}$
is the most accurate in both norms at all four targets. Using
$h_{\mathrm{eff}}\propto N_T^{-1/2}$ and the actual coarsest and finest
mesh counts, the observed $L^2$ and $H^1$-seminorm convergence rates
are approximately $0.88$ and $0.69$ for $\mathbb M_H$, compared with
$1.14$ and $1.06$ for $\mathbb M_{H+\mathrm{adapt}}$.
The $L^2$ rates are below the nominal second-order spatial rate for
$P_1$ elements. Because the computations use a fixed time step and
repeated remeshing with solution transfer, the observed rates reflect
the complete numerical procedure rather than spatial discretization
error alone.

We focus here on $\mathbb M_H$ and $\mathbb M_{H+\mathrm{adapt}}$
to assess the effect of adding a solution-dependent metric density.
In two dimensions the pure effective-tensor metric has unit determinant
at the element level, whereas $\mathbb M_{H+\mathrm{adapt}}$ multiplies
it by the recovered-Hessian factor in~\eqref{M-H-adapt-2d}.
At target $N\approx2048$, the $L^2$ errors are $7.625\times10^{-3}$ and $5.307\times10^{-3}$, respectively,
and the $H^1$-seminorm errors are $0.2523$ and $0.1881$.
At target $N\approx4608$, the corresponding errors are
$5.626\times10^{-3}$ and $3.678\times10^{-3}$ in $L^2$, and
$0.1850$ and $0.1349$ in the $H^1$ seminorm.
Thus the combined metric reduces both errors relative to the pure
$\mathbb M_H$ at these budgets, though it does not match the accuracy
of Hessian-only adaptation.

\begin{figure}[htbp]
\centering
\subfloat[Pure $\mathbb M_H$, $N_T=798$]{%
  \manuscriptfigure{.47\textwidth}{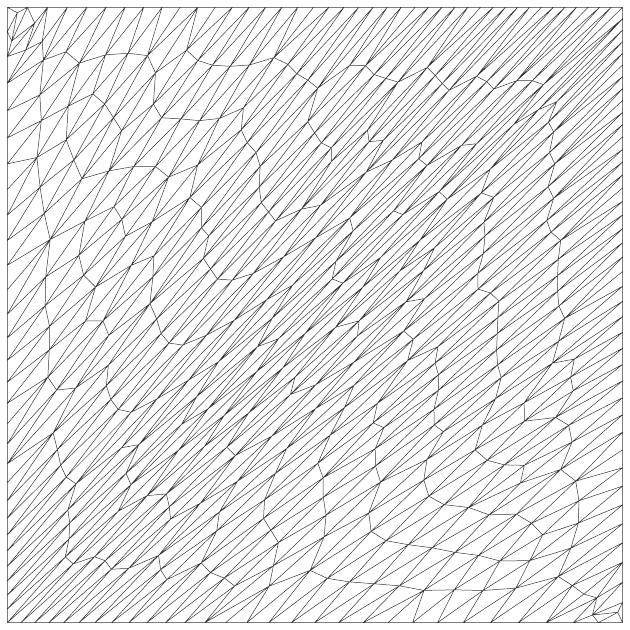}}\hfill
\subfloat[$\mathbb M_{H+\mathrm{adapt}}$, $N_T=874$]{%
  \manuscriptfigure{.47\textwidth}{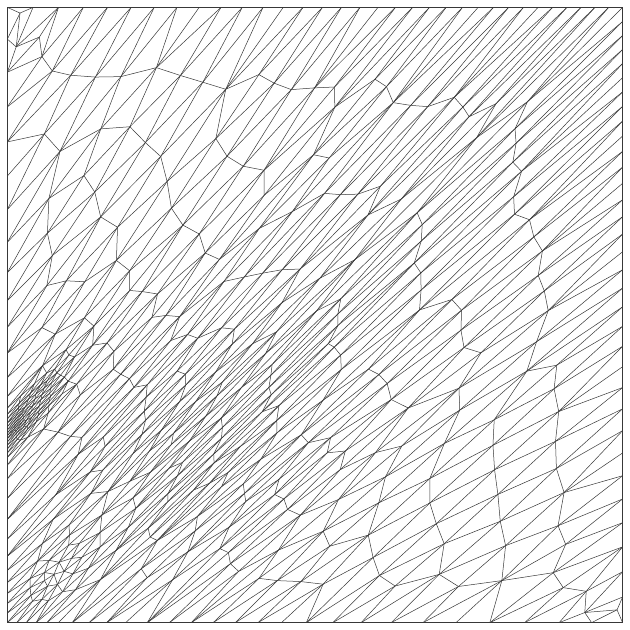}}
\caption{Example~3: accepted final meshes at $T=0.3$ and nominal target
$N\approx800$. The determinant-normalized $\mathbb M_H$ uses the
SUPG-effective-tensor geometry; $\mathbb M_{H+\mathrm{adapt}}$ adds
solution-dependent density through the recovered Hessian.
The actual final triangle counts differ.}
\label{fig:ex3-density-meshes}
\end{figure}

The layer is steepest in the direction $(1,1)^T/\sqrt2$,
an eigenvector of $\mathbb D$ associated with eigenvalue $100$.
The coefficient-dependent orientation alone does not necessarily
resolve the layer efficiently. The lower errors for
$\mathbb M_{H+\mathrm{adapt}}$ illustrate the benefit of combining
effective-tensor geometry with solution-dependent density for this
problem. Since the equation has a manufactured source, this example
assesses the role of metric density in approximation accuracy, not
DMP satisfaction.

\subsection{Numerical verification of sufficient DMP conditions}
\label{subsec:fixed-certification}
The adaptive transport tests measure numerical bound violations, while
Theorem~\ref{thm-DMP} provides sufficient conditions for their absence.
To illustrate that these conditions are attainable, we consider an
independent fixed-mesh problem on a transformed rhombus. Let
\[
\mathbb D=\begin{pmatrix}50.5&49.5\\49.5&50.5\end{pmatrix},\qquad
\Omega=J(0,1)^2,\qquad
J=\frac{1}{\sqrt{10}}\mathbb D^{1/2}
\begin{pmatrix}1&-1/2\\0&\sqrt3/2\end{pmatrix}.
\]
Writing $(s,r)^T=J^{-1}\boldsymbol x$, we set
$\psi(s,r)=\sin^2(\pi s)\sin^2(\pi r)$ and
\[
\boldsymbol b(\boldsymbol x)=\frac{1}{\pi\sqrt2\|J\|_2}
J\begin{pmatrix}\partial_r\psi\\-\partial_s\psi\end{pmatrix}.
\]
The velocity is divergence free and vanishes on the boundary. The
nonnegative initial hill is
$u_0(\boldsymbol x)=\tfrac12[1+\cos(4\pi\rho)]$ when
$\rho=\|J^{-1}\boldsymbol x-(1/2,1/2)^T\|_2\le1/4$, and zero
otherwise. Both the source and the Dirichlet data vanish. We use
backward Euler with limited SUPG stabilization.

\begin{table}[htbp]
\centering\small
\begin{tabular}{rrrrrr}
\toprule
$\varepsilon$ & Triangles & $\Delta t$ & $\delta_{\min}$ & $\Delta t_L$ & Observed $(E_-,E_+)$ \\
\midrule
$0.003$ & $5000$ & $0.010417$ & $-8.864$ & -- & $0,\ 0$ \\
$0.003$ & $9800$ & $0.010417$ & $15.591$ & $0.009228$ & $0,\ 0$ \\
$0.010$ & $5000$ & $0.005$ & $107.803$ & $0.001350$ & $0,\ 0$ \\
$0.010$ & $9800$ & $0.005$ & $244.257$ & $0.000626$ & $0,\ 0$ \\
\bottomrule
\end{tabular}
\caption{Independent fixed-mesh DMP checks on the transformed rhombus.
A negative $\delta_{\min}$ indicates failure of the sufficient spatial
condition, in which case no sufficient lower time-step bound is
certified.}
\label{tab:fixed-dmp-certification}
\end{table}

Table~\ref{tab:fixed-dmp-certification} shows that for
$\varepsilon=0.003$ and $9800$ triangles, the mass, spatial,
assembled-sign, and time-step requirements all hold, and no nodal
bound violation is observed. The corresponding $5000$-triangle mesh
fails the sufficient spatial condition, although its computed solution
remains within the bounds. This illustrates that failure of a
sufficient condition does not imply a violation for a particular
trajectory. This fixed-mesh test verifies that the sufficient
conditions are jointly attainable on a nontrivial anisotropic mesh.
For changing meshes, the corresponding assembled-matrix and time-step
conditions, together with the bound-preserving transfer requirement,
must hold at each step. The BAMG-adaptive transport experiments below
instead assess observed bound violations without claiming that their
complete trajectories satisfy these sufficient conditions.

\subsection{Example 4: rotating-pulse transport with moderate diffusion anisotropy}
\label{subsec:e4-supg-metric}
Finally, we study bound violations with anisotropic physical diffusion and use a 
controlled comparison to examine the influence of the effective tensor.
On $\Omega=(-1,1)^2$, consider
\begin{equation}\label{eq:e4-pde}
u_t-\varepsilon\nabla\cdot(\mathbb D\nabla u)+\boldsymbol b\cdot\nabla u=0,
\quad\varepsilon=10^{-3},\quad
\mathbb D=\begin{pmatrix}3&2\\2&3\end{pmatrix},\quad
\boldsymbol b=6\begin{pmatrix}y(1-x^2)\\-x(1-y^2)\end{pmatrix}.
\end{equation}
The constant SPD diffusion tensor has eigenvalues $5$ and $1$ and
principal axes at $45^\circ$ to the coordinate axes. The velocity is
divergence free and boundary tangent. Its streamlines follow the
level sets of $\psi(x,y)=\frac12(1-x^2)(1-y^2)$. Precisely,
$\boldsymbol b=6(-\partial_y\psi,\partial_x\psi)^T$.
Homogeneous Dirichlet conditions
are imposed. The compact interior initial pulse is
\begin{equation}\label{eq:e4-pulse}
C(s)=\begin{cases}\frac12(1+\cos(\pi s)),&|s|<1,\\0,&|s|\ge1,\end{cases}
\qquad
u_0(x,y)=C\!\left(\frac{\psi(x,y)-0.34875}{0.045}\right)
 C\!\left(\frac{\operatorname{atan2}(y,x)}{0.75}\right),
\end{equation}
with $u_0(0,0)=0$. Thus $0\le u_0\le1$ and the initial support is
strictly interior. The effective tensor used by the mesh metric is
$\mathbb H_K=\varepsilon\mathbb D+
\tau_K\overline{\boldsymbol b\boldsymbol b^T}_K$, with the actual
limited stabilization parameter
$\tau_K=\min\{\tau_K^{\mathrm{nominal}},
[4B_K\max_i\|\nabla\varphi_i\|]^{-1}\}$ used by the PDE operator.
The diffusion-based metric employs this anisotropic $\mathbb D$, not $I$.

All principal runs use backward Euler with $\Delta t=0.003125$
to $T=0.25$, with two prescribed mesh updates at each scheduled
physical remeshing time. The five methods share BAMG settings,
limited SUPG stabilization and the stated transfer policy.
The results are presented in Table \ref{tab:e4-moderate-anisotropy}.

\begin{table}[htbp]
\centering\scriptsize\setlength{\tabcolsep}{2.5pt}
\begin{tabular}{lrrrrrrrr}
\toprule
Metric & $N_T$ & $E_-^{\max}$ & $I_-$ & $I_{-,2}$ & $E_+^{\max}$ & $L^2_{\mathrm{ref}}(\%)$ & $H^1_{\mathrm{ref}}(\%)$ & Peak/ref $(\%)$\\
\midrule
\multicolumn{9}{l}{Target 2048 triangles}\\
FM--SUPG & 2048 & 0.064191 & 0.008688 & 0.00064841 & 0 & 26.176 & 54.439 & 67.443\\
$\mathbb{M}_{\mathrm{adapt}}$ & 2082 & 0.0030994 & 0.00031140 & 2.5018e-05 & 0 & 12.377 & 15.838 & 86.501\\
$\mathbb{M}_{\mathrm{DMP+adapt}}$ & 2058 & 0.0044878 & 0.00055819 & 4.3100e-05 & 0 & 16.413 & 24.656 & 82.073\\
${\mathbb{M}_H}$ & 2048 & 0.0065249 & 0.0011262 & 7.3129e-05 & 0 & 20.899 & 38.83 & 75.044\\
$\mathbb{M}_{H+\mathrm{adapt}}$ & 2040 & 0.0019962 & 0.00012673 & 9.1458e-06 & 0 & 15.003 & 20.599 & 83.913\\
\addlinespace[2pt]
\multicolumn{9}{l}{Target 5000 triangles}\\
FM--SUPG & 5000 & 0.012806 & 0.0017053 & 0.00011239 & 0 & 16.203 & 34.834 & 79.918\\
$\mathbb{M}_{\mathrm{adapt}}$ & 5021 & 0.0016764 & 0.00012372 & 7.3795e-06 & 0 & 10.856 & 12.926 & 88.713\\
$\mathbb{M}_{\mathrm{DMP+adapt}}$ & 5004 & 0.0023575 & 0.00022952 & 1.3409e-05 & 0 & 12.273 & 16.777 & 86.954\\
${\mathbb{M}_H}$ & 5054 & 0.00459713 & 0.00025527 & 9.7375e-06 & 0 & 21.767 & 37.363 & 75.89\\
$\mathbb{M}_{H+\mathrm{adapt}}$ & 5000 & 0.00073741 & 5.4570e-05 & 3.2150e-06 & 0 & 11.952 & 15.046 & 87.586\\
\bottomrule\end{tabular}
\caption{Example~4, $\mathbb D=\left(\begin{smallmatrix}3&2\\2&3\end{smallmatrix}\right)$,
$\varepsilon=10^{-3}$, $\Delta t=0.003125$ and limited SUPG.
Accepted-time DMP maxima and integrals, relative independent numerical-
reference errors ($H^1$ gradient seminorm), and final peak retention.
All five methods use the paired principal protocol. Pure $\mathbb{M}_H$ is determinant-normalized in two dimensions. 
Adaptive methods use sign-aware, globally mass-preserving transfer.
$N_T$ is the actual final triangle count.}
\label{tab:e4-moderate-anisotropy}
\end{table}

The combined metric $\mathbb{M}_{H+\mathrm{adapt}}$ reduces $E_-^{\max}$ relative to normalized
$\mathbb{M}_H$ by $69.41\%$ and $83.96\%$ at the two budgets.
Neither method satisfies all sufficient assembled-matrix DMP conditions.
Relative to $\mathbb{M}_{\mathrm{DMP+adapt}}$, the combined metric
reduces maximum undershoot by $55.5\%$ and $68.7\%$ at the two
budgets. Relative to $\mathbb{M}_{\mathrm{adapt}}$ the reductions
are $35.6\%$ and $56.0\%$. The corresponding time-integrated
violations also decrease. The metric ranking and its change with mesh budget are problem dependent.

\begin{figure}[htbp]
\centering
\subfloat[$\mathbb M_{\mathrm{adapt}}$ mesh, $N_T=5021$]{%
  \manuscriptfigure{.47\textwidth}{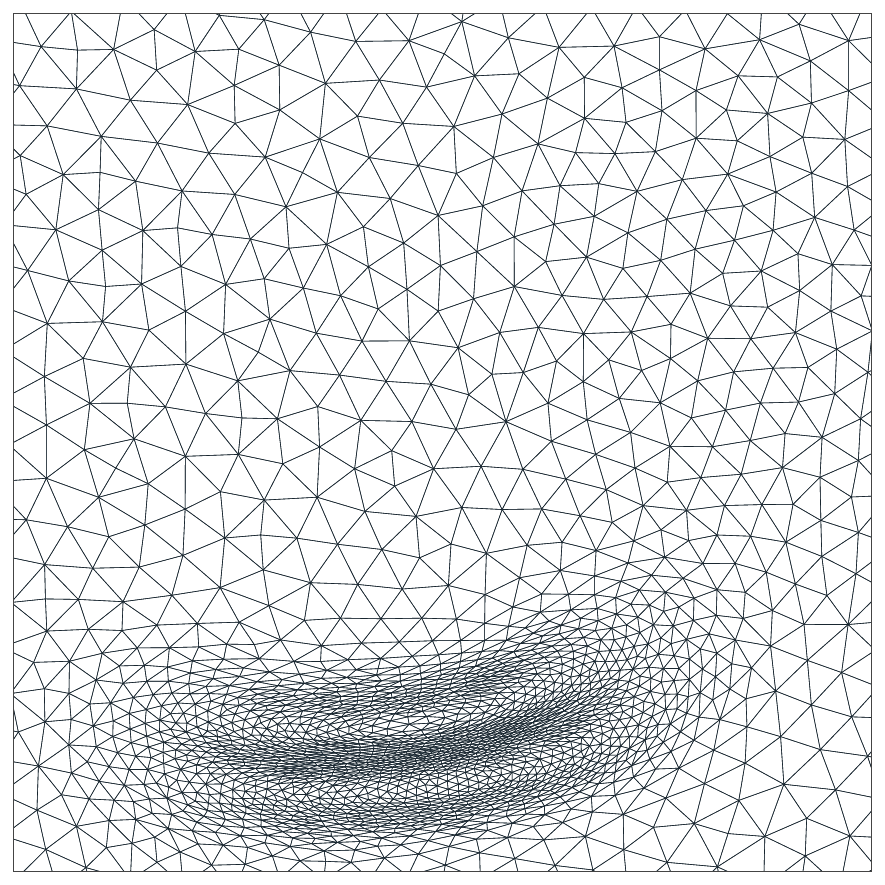}}
\subfloat[$\mathbb M_{H+\mathrm{adapt}}$ mesh, $N_T=5000$]{%
  \manuscriptfigure{.47\textwidth}{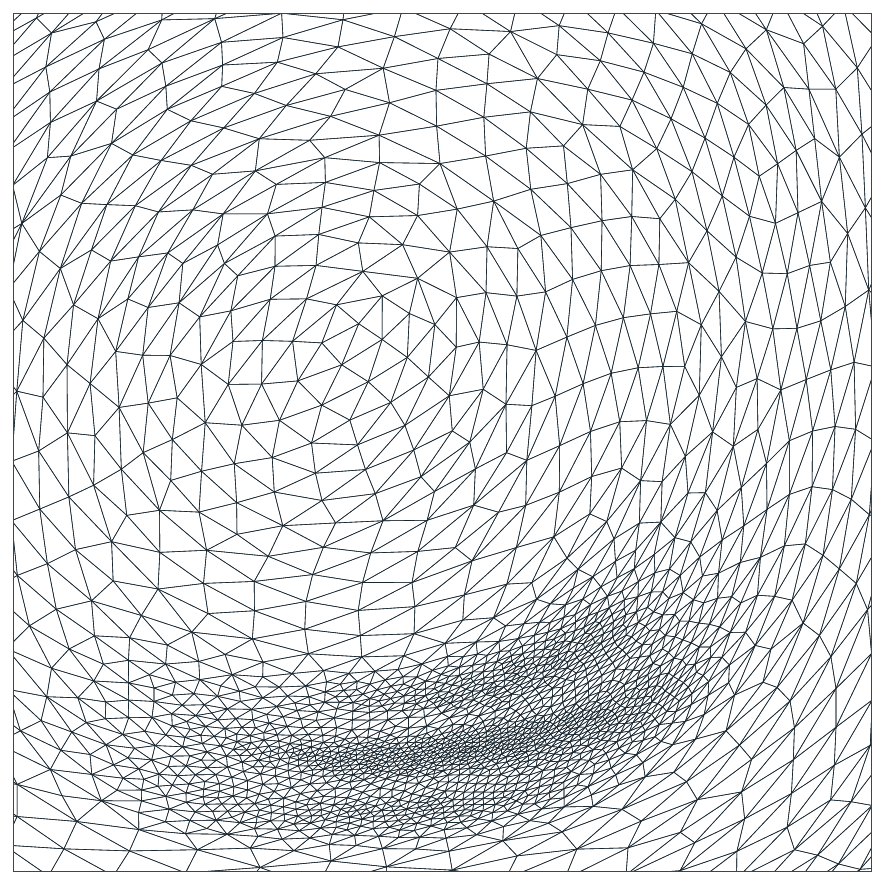}}\\[2pt]
\subfloat[$\mathbb M_{\mathrm{adapt}}$ solution]{%
  \manuscriptfigure{.47\textwidth}{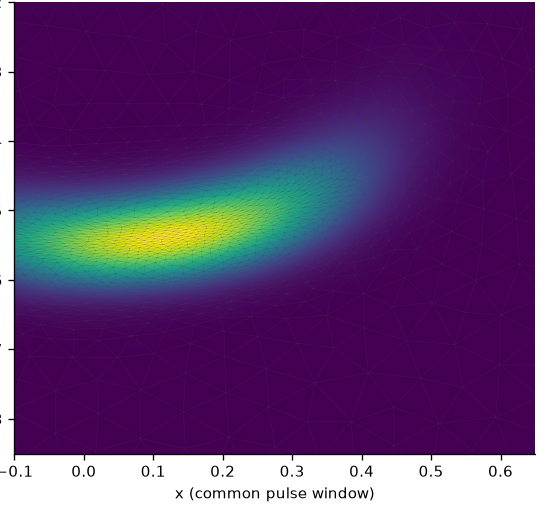}}
\subfloat[$\mathbb M_{H+\mathrm{adapt}}$ solution]{%
  \manuscriptfigure{.47\textwidth}{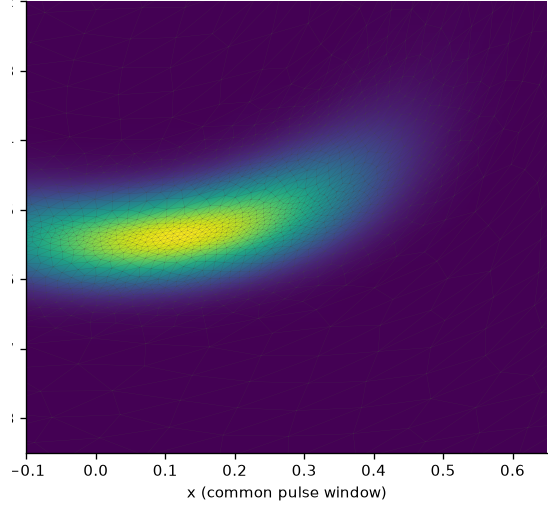}}
\caption{Example~4 at $T=0.25$ and target $N\approx5000$. The top
row shows accepted adapted meshes for the two metrics with the same
zoom and coordinate scale. The bottom row shows the corresponding
accepted finite element solutions using a common color scale.
The mesh panels use $-0.45\le x\le0.75$, $-0.7\le y\le0.5$;
the solution panels show a tighter view of the pulse.}
\label{fig:e4-mesh-solution-combined}
\end{figure}

At the finer budget, $\mathbb{M}_{H+\mathrm{adapt}}$ has relative
reference errors of $11.95\%$ ($L^2$) and $15.05\%$ ($H^1$),
compared with $10.86\%$ and $12.93\%$ for Hessian-only adaptation.
The independent Rannacher--Crank--Nicolson reference uses the same
diffusion tensor, a finest mesh of $3{,}276{,}800$ triangles, and
$\Delta t=0.00078125$. The last spatial changes were $0.0387\%$
in $L^2$ and $1.89\%$ in the $H^1$ seminorm. The last temporal
changes were below $0.05\%$ in both norms. 
The meshes and corresponding solutions are compared in
Figure~\ref{fig:e4-mesh-solution-combined}.

To distinguish the SUPG-effective tensor substitution from the
independently prescribed mesh density, we performed a four-cell
comparison at the two Example~4 complexities. For
$\mathbb A_K\in\{\mathbb D_K,\mathbb H_K\}$ define in two dimensions
\begin{equation}\label{eq:tensor-density-controlled}
\mathbb{M}_{A,\rho,K}=C\rho_K\det(\mathbb A_K)^{1/2}\mathbb A_K^{-1},
\end{equation}
where $\rho$ is either constant or the same externally fixed,
Hessian-derived physical-space density field for both tensor choices.
The common density is formed once from a prior independent reference
trajectory, rather than recalculated for each candidate adaptive
solution. Thus this test does not identify the published
$\mathbb{M}_{H+\mathrm{adapt}}$ with the H/common cell.
The eight trajectories use the same physical remeshing schedule and
limited-SUPG/sign-aware protocol. The results are given in Table~\ref{tab:tensor-density-controlled}.

\begin{table}[htbp]
\centering\small
\begin{tabular}{lrr}
\toprule
Tensor/density & $E_-^{\max}$, $N\approx2048$ & $E_-^{\max}$, $N\approx5000$\\
\midrule
D/uniform & 0.060196 & 0.012702\\
D/common & 0.023770 & 0.011513\\
H/uniform & 0.006525 & 0.003674\\
H/common & 0.003581 & 0.001565\\
\bottomrule
\end{tabular}
\caption{Example~4 controlled effective-tensor and prescribed-density
comparison.}
\label{tab:tensor-density-controlled}
\end{table}

At common density, the replacement $\mathbb D\to\mathbb H$
reduces $E_-^{\max}$ from $0.0237696$ to $0.00358059$
($84.94\%$) at target $2048$, and from $0.0115126$ to
$0.00156473$ ($86.41\%$) at target $5000$. These observations
support a contribution from SUPG-effective tensor geometry in this
controlled deterministic comparison. However, $\mathbb H$ also
changes the local SUPG limiter, and the H/common cases have more
forbidden positive off-diagonal $B$ entries than the D/common cases
($6030$ versus $5664$, and $14577$ versus $13311$, respectively).
Improved undershoot amplitude is therefore not evidence that
this substitution approaches assembled-matrix DMP certification.
The single-trajectory density interactions are descriptive and are not
interpreted as robust or statistically resolved effects. No universal tensor
ranking is inferred.

We also evaluated mass and assembled-step matrix signs throughout the
accepted transport trajectories. All tested assembled operators
retain forbidden positive off-diagonal entries. Hence none of the
adaptive trajectories is certified by our sufficient DMP theorem.
Mass-compatible limited stabilization substantially reduces the
nominal streamline parameter on many elements. 
This contrasts with the independent fixed-mesh certification in
Section~\ref{subsec:fixed-certification}.

Fully adaptive time-step refinements at the same physical remeshing
times preserve the combined metric's smallest undershoot among the five
methods at both budgets for Examples~2 and~4 and all three tested
steps. The narrowest comparison is Example~4 at target $N=5000$ and
$\Delta t/4$: the combined metric gives $0.00150312$ against
$0.00155816$ for Hessian-only adaptation ($3.53\%$ lower), with
actual final counts of $4900$ and $5022$ triangles, respectively.
Thus the ranking is descriptive, not a mesh-count-independent or
asymptotic convergence conclusion. 

\section{Conclusions}
\label{sec:conclusions}
We derived sufficient fully discrete nodal DMP conditions for a
$P_1$ SUPG method on anisotropic meshes. The analysis accounts for the
stabilized mass matrix, the effective tensor
$\mathbb H_K=\varepsilon\mathbb D_K+\tau_K\overline{\boldsymbol b\boldsymbol b^T}_K$,
time-step restrictions, and bound-preserving transfer after remeshing.
An assembled-matrix criterion complements the conservative local
geometric conditions.

The resulting determinant-normalized metric $\mathbb M_H$ uses the
SUPG-effective-tensor geometry. The combined metric
$\mathbb M_{H+\mathrm{adapt}}$ retains that geometry while introducing
solution-dependent mesh density. These constructions aim to improve
the behavior of a linear SUPG discretization through anisotropic mesh
adaptation, rather than through an additional nonlinear
bound-enforcing modification of the PDE operator.

In the tested isotropic and moderately anisotropic transport problems,
$\mathbb M_{H+\mathrm{adapt}}$ reduces maximum-in-time undershoots
relative to the comparison metrics under mass-compatible stabilization.
The controlled anisotropic experiment indicates a contribution from
the effective tensor when a common density is prescribed. The smooth
manufactured examples also show that Hessian-only adaptation can be
more accurate, and the transport calculations reveal an approximation
accuracy cost for some DMP-oriented meshes.

The reported advantages depend on the numerical problem, stabilization
policy, and time-stepping protocol. The adaptive transport computations
do not satisfy all sufficient algebraic DMP conditions, so smaller
violations must not be interpreted as a guarantee of bound preservation.
The derived conditions nevertheless identify geometric and algebraic
requirements relevant to DMP-oriented mesh design and provide a basis
for further improvements in stabilization and adaptation.

\paragraph{Data availability.}
Numerical data supporting the tables and figures accompany the
supplementary material. Complete raw simulation trajectories are not
publicly deposited.

\paragraph{Declaration of AI-assisted tools.}
AI-assisted tools (ChatGPT and Codex) were used in manuscript editing,
code-related tasks, and preparation of numerical reporting materials.
The authors take responsibility for the mathematical analysis,
verification and interpretation of the numerical evidence, and final
manuscript content.

\end{document}